\documentclass[12pt]{amsart}

\usepackage{amsmath}
\usepackage{amscd}
\usepackage{amssymb}
\usepackage{latexsym}
\usepackage{stmaryrd}
\usepackage{color}

\usepackage[arrow,curve,matrix,tips,2cell]{xy}

\usepackage[hmargin=3cm,vmargin=2.5cm,includehead,includefoot]{geometry}

\usepackage{fancyhdr}

\newtheorem{theorem}{Theorem}[section]
\newtheorem*{theorem*}{Theorem}
\newtheorem{lemma}[theorem]{Lemma}
\newtheorem*{lemma*}{Lemma}
\newtheorem{corollary}[theorem]{Corollary}
\newtheorem*{corollary*}{Corollary}
\newtheorem{proposition}[theorem]{Proposition}

\newtheorem{remark}[theorem]{Remark}
\newtheorem{question}[theorem]{Question}
\newtheorem{definition}[theorem]{Definition}

\newtheorem{example}[theorem]{Example}
\newtheorem{examples}[theorem]{Examples}
\newcommand{\bgl}{\begin{equation}} %eine Gleichung mit Ziffer
\newcommand{\egl}{\end{equation}}
\newcommand{\bgloz}{\begin{equation*}} %eine Gleichung ohne Ziffer
\newcommand{\egloz}{\end{equation*}}
\newcommand{\bgln}{\begin{eqnarray}} %mehrere Gleichungen mit Ziffer
\newcommand{\egln}{\end{eqnarray}}
\newcommand{\bglnoz}{\begin{eqnarray*}} %mehrere Gleichungen ohne Ziffer
\newcommand{\eglnoz}{\end{eqnarray*}}
\newcommand{\btheo}{\begin{theorem}}
\newcommand{\etheo}{\end{theorem}}
\newcommand{\btheooz}{\begin{theorem*}}
\newcommand{\etheooz}{\end{theorem*}}
\newcommand{\blemma}{\begin{lemma}}
\newcommand{\elemma}{\end{lemma}}
\newcommand{\blemmaoz}{\begin{lemma*}}
\newcommand{\elemmaoz}{\end{lemma*}}
\newcommand{\bproof}{\begin{proof}}
\newcommand{\eproof}{\end{proof}}
\newcommand{\bbew}{\begin{beweis}}
\newcommand{\ebew}{\end{beweis}}
\newcommand{\bremark}{\begin{remark}\em}
\newcommand{\eremark}{\end{remark}}
\newcommand{\bquestion}{\begin{question}\em}
\newcommand{\equestion}{\end{question}}
\newcommand{\bdefin}{\begin{definition}}
\newcommand{\edefin}{\end{definition}}
\newcommand{\bprop}{\begin{proposition}}
\newcommand{\eprop}{\end{proposition}}
\newcommand{\bcor}{\begin{corollary}}
\newcommand{\ecor}{\end{corollary}}
\newcommand{\bcoroz}{\begin{corollary*}}
\newcommand{\ecoroz}{\end{corollary*}}
\newcommand{\bfa}{\begin{cases}} %Fallunterscheidung
\newcommand{\efa}{\end{cases}}
\newcommand{\bexample}{\begin{example}\em}
\newcommand{\eexample}{\end{example}}
\newcommand{\bexamples}{\begin{examples}\em}
\newcommand{\eexamples}{\end{examples}}
\newcommand{\cA}{\mathcal A}
\newcommand{\cB}{\mathcal B}
\newcommand{\cC}{\mathcal C}
\newcommand{\cD}{\mathcal D}

\newcommand{\cH}{\mathcal H}

\newcommand{\cK}{\mathcal K}
\newcommand{\cL}{\mathcal L}

\newcommand{\cN}{\mathcal N}
\newcommand{\cO}{\mathcal O}

\newcommand{\cS}{\mathcal S}

\newcommand{\cV}{\mathcal V}

\newcommand{\cZ}{\mathcal Z}
\def\Cz{\mathbb{C}}
\def\Fz{\mathbb{F}}

\def\Nz{\mathbb{N}}

\def\Rz{\mathbb{R}}
\def\Tz{\mathbb{T}}
\def\Zz{\mathbb{Z}}

\def\1z{\mathbb{1}}

\newcommand{\bA}{\mathbf{A}}

\newcommand{\bM}{\mathbf{M}}
\newcommand{\an}[1]{``#1''} % Anfuehrungsstriche
\newcommand{\ti}{\tilde}

\newcommand{\lori}{\longrightarrow}

\newcommand{\ma}{\mapsto} % wird abgebildet auf
\newcommand{\mafr}{\mapsfrom} % kommt von
\newcommand\onto{\twoheadrightarrow} % surjektiv
\newcommand\into{\hookrightarrow} % injektiv
\newcommand{\Rarr}{\Rightarrow} % Folgerung
\newcommand{\LRarr}{\Leftrightarrow} % aequivalent

\newcommand{\ve}{\varepsilon}

\def\SEMI{\mbox{$\times\kern-2pt\vrule height5pt width.6pt \kern3pt $}}

\newcommand{\Aut}{{\rm Aut}\,}

\newcommand{\img}{{\rm im\,}}

\newcommand{\Spec}{{\rm Spec\,}} % Spektrum
\newcommand{\rk}{{\rm rk\,}}

\newcommand{\id}{{\rm id}}

\renewcommand{\dim}{{\rm dim}\,}

\newcommand{\abs}[1]{\left|#1\right|} % Betrag
\newcommand{\norm}[1]{\left\|#1\right\|} % Norm
\newcommand{\defeq}{\mathrel{:=}} % per Definition
\newcommand{\dop}{\text{: }} % in Mengen
\newcommand{\supp}{{\rm supp}\,}

\newcommand{\lge}{\left\{} % links geschweift
\newcommand{\rge}{\right\}} % rechts geschweift
\newcommand{\lru}{\left(} % links rund
\newcommand{\rru}{\right)} % rechts rund
\newcommand{\lsp}{\left\langle} % links spitz
\newcommand{\rsp}{\right\rangle} % links spitz
\newcommand{\rukl}[1]{\lru #1 \rru} % runde Klammer
\newcommand{\gekl}[1]{\lge #1 \rge} % geschweifte Klammer
\newcommand{\spkl}[1]{\lsp #1 \rsp} % spitze Klammer

\newcommand{\menge}[2]{\gekl{ #1 \, \dop #2 }} % Menge
\title[Cartan subalgebras in C*-algebras. Existence and uniqueness]{Cartan subalgebras in C*-algebras. \\ Existence and uniqueness}

\author{Xin Li}
\address{School of Mathematical Sciences\\
Queen Mary University of London\\
Mile End Road\\
London E1 4NS\\
United Kingdom}
\email{xin.li@qmul.ac.uk}

\author{Jean Renault}
\address{Universit\'e d'Orl\'eans et CNRS  (UMR 7349 et FR2964), D\'epartement de Math\'ematiques, 
F-45067 Orl\'eans Cedex 2, France}
\email{jean.renault@univ-orleans.fr}


\thanks{The first named author is supported by EPSRC grant EP/M009718/1.}

\begin{document}

\begin{abstract}
We initiate the study of Cartan subalgebras in C*-algebras, with a particular focus on existence and uniqueness questions. For homogeneous C*-algebras, these questions can be analysed systematically using the theory of fibre bundles. For group C*-algebras, while we are able to find Cartan subalgebras in C*-algebras of many connected Lie groups, there are classes of (discrete) groups, for instance non-abelian free groups, whose reduced group C*-algebras do not have any Cartan subalgebras. Moreover, we show that uniqueness of Cartan subalgebras usually fails for classifiable C*-algebras. However, distinguished Cartan subalgebras exist in some cases, for instance in nuclear Roe algebras.
\end{abstract}

\maketitle

%\tableofcontents

\setlength{\parindent}{0cm} \setlength{\parskip}{0.25cm}

\section{Introduction}

The construction of groupoid C*-algebras is an extremely powerful and general method to produce C*-algebras. Many C*-algebras admit descriptions as groupoid C*-algebras, either directly by construction or indirectly by -- often explicitly given -- identifications. Having groupoid models at hand is very helpful to analyse the C*-algebraic structure, be it for determining ideals, for nuclearity, computing K-theory, the UCT etc. With this in mind, it is natural to ask the following questions: Which C*-algebras admit groupoid models? And to what extent are these groupoid models unique? As we shall explain below, these two questions are closely related to the existence and uniqueness questions for Cartan subalgebras in C*-algebras.

Our goal in this paper is to initiate a systematic study of Cartan subalgebras in C*-algebras, with a particular focus on existence and uniqueness. Classification of C*-algebras has seen tremendous advances recently \cite{W10, W12, MS12, MS14, SWW, GLN, EGLN, TWW}. However, finer structures given by Cartan subalgebras have not yet been analysed in detail. Apart from being of interest for C*-algebra theory itself, the notion of Cartan subalgebras also builds a bridge to topological dynamical systems, as it is closely related to the notion of continuous orbit equivalence \cite{Li_COER}. The latter, in turn, establishes a connection to geometric group theory via dynamic characterizations of quasi-isometry \cite{Li_DQH}. The loop is closed by the observation in Theorem~\ref{Roe-distCartan}, building on \cite{SW,Li_DQH}, that quasi-isometry can also be characterized using Cartan subalgebras in Roe algebras.

The situation we are describing here is completely analogous to the measurable framework, where a fruitful interplay between Cartan subalgebras in von Neumann algebras (vN-algebras), orbit equivalence for measure-preserving dynamical systems and geometric group theory has been established \cite{Sha05,Gab10,Fur11}. Indeed, the comparison between C*-algebraic and vN-algebraic Cartan subalgebras is one of the leading principles behind our present work.

Coming back to our original questions on existence and uniqueness of groupoid models, we have to make our question more precise in order to get a satisfactory answer.
%added:
%----------------------------------------
%The following question is often asked when one meets the notion of groupoid C*-algebras: is every C*-algebra a groupoid C*-algebra? This has to be made more precise to get a satisfactory answer. 
What kind of groupoids should be allowed? For example, group C*-algebras, which are amongst the most challenging C*-algebras, are groupoid C*-algebras by definition; however in some cases, e.g. abelian or $ax+b$-groups, they may admit more revealing groupoid realisations. Usually, one looks for a groupoid which has as little isotropy as possible. In the case of vN-algebras, it suffices to consider groupoids without isotropy (also known as principal groupoids or equivalence relations). In the theory of C*-algebras, we allow some isotropy and consider topologically principal groupoids (this means that the points without isotropy are dense) to cover C*-algebras arising from foliations or from pseudo-groups. In this paper, we work with groupoids which are topologically principal, locally compact, Hausdorff and \'etale. The reason why we choose such groupoids is that the reduced C*-algebras of such groupoids, equipped possibly with a twist, admit a convenient C*-algebraic characterisation, namely the existence of a Cartan subalgebra, a notion borrowed from the theory of vN-algebras: An abelian sub-C*-algebra $B$ of a C*-algebra $A$ is called a {\it Cartan subalgebra} if it contains an approximate unit of $A$, is maximal abelian, regular and there exists a faithful conditional expectation of $A$ onto $B$. By \cite{R08}, every Cartan pair $(A,B)$ consisting of a Cartan subalgebra $B$ in a C*-algebra $A$ is of the form $(C^*_r(G,\Sigma), C_0(G^{(0)}))$, where the groupoid $G$ is as above, $G^{(0)}$ is its unit space and $\Sigma$ is a twist over $G$. The case when $G$ is an equivalence relation had been worked out earlier by A. Kumjian: $G$ has no isotropy if and only if the subalgebra $B$ has in addition the unique extension property; then it is called a {\it diagonal}. The original definition of a Cartan subalgebra in the theory of vN-algebras is an abelian subalgebra $\bA$ of a vN-algebra $\bM$ which is a masa, regular and there exists a faithful normal contional expectation of $\bM$ onto $\bA$. The generic Cartan pair is $({\bM},{\bA})=(\bM(R,\sigma), L^\infty(X))$, where $R$ is a countable standard measured equivalence $R$ on $X$ and $\sigma$ is a 2-cocycle (see \cite{FMII}). Existence and uniqueness of Cartan subalgebras in vN-algebras is an active field of research \cite{Voi,Oza,OP10a,OP10b,PV14a,PV14b,BHV}. The examples presented in this paper will show that the investigation of Cartan subalgebras in C*-algebras is no less interesting. It turns out that, unless a C*-algebra contains a Cartan subalgebra from its very construction, showing the existence of a Cartan subalgebra involves a variety of techniques. It will also turn out that, when Cartan subalgebras exist, there are usually infinitely many of them, at least in case the big C*-algebra is classifiable. Let us briefly describe the content of the paper section by section.

The second section studies Cartan subalgebras in $n$-homogeneous C*-algebras. When $n$ is finite, it is shown that a Cartan pair is necessarily locally trivial. This puts the problem into the general theory of fibre bundles: the $n$-homogeneous C*-algebra $A$ admits a Cartan subalgebra if and only the principal ${\rm Aut}(M_n)$-bundle associated to $A$ can be reduced to the subgroup ${\rm Aut}(M_n,D_n)$. Using this theory, we give a complete enumeration of Cartan subalgebras of an $n$-homogeneous C*-algebra over a $k$-sphere. This gives a conceptual explanation of the result of  K.D. Gregson \cite{Gre} and T. Natsume (see appendix of \cite{Kum}) that for $k>2$, only the trivial $n$-homogeneous C*-algebra $C(S^k, M_n)$ admits a Cartan subalgebra; moreover, our study shows that this Cartan subalgebra is essentially unique. A similar study is made when $n$ is infinite, under the assumption that the pair $(A,B)$ is locally trivial.

The third section gives the existence of Cartan subalgebras in group C*-algebras of certain connected Lie groups, building on \cite{Fel,Mil,KM,Val,PP}. Namely, it is shown that the full group C*-algebras of $SL(2,\Cz)$, $SL(2,\Rz)$ as well as the reduced group C*-algebras of every connected semi-simple Lie group with real rank one and finite center and of every connected complex semi-simple Lie group have a Cartan subalgebra. The proof uses the description of these C*-algebras as sub-C*-algebras of $C_0(M, {\mathcal K}({\mathcal H}))$, where $M$ is a locally compact Hausdorff space and $\mathcal H$ is an infinite dimensional Hilbert space.

In the fourth section, we show that certain reduced group C*-algebras do not have Cartan subalgebras. This is in particular the case for non-abelian free groups. Here the strategy is to use the known corresponding results for the group vN-algebras \cite{Voi,Hay,Oza,OP10a,OP10b,PV14a,PV14b,BHV}. Note that some extra work is needed to show the non-existence of Cartan subalgebras rather than diagonals.

The fifth section illustrates various techniques to construct infinitely many Cartan subalgebras in a given C*-algebra. The first one uses the recent construction of Deeley, Putnam and Strung \cite{DPS} of the Jiang-Su algebra $\mathcal Z$ as a groupoid C*-algebra. It applies to unital C*-algebras which have a Cartan subalgebra whose spectrum has finite covering dimension and which are $\mathcal Z$-stable. The second one, described by A. Kumjian in \cite{Kum84}, uses non-stable K-theory and applies to irrational rotation C*-algebras. In both cases, the Cartan subalgebras are shown to be non-isomorphic. The last result is that unital (resp. stable) UCT Kirchberg algebras have infinitely many inequivalent Cartan subalgebras which are all isomorphic to the algebra of continuous functions on the Cantor space (or its locally compact non-compact analogue). It relies on a groupoid realisation of these algebras due to J. Spielberg \cite{Sp07a}.

The phenomena presented in Section~5 show that in general, we cannot expect Cartan subalgebras in C*-algebras to be rigid, i.e., uniqueness fails. However, it turns out that in some cases, we can find distinguished Cartan subalgebras. This means that within a class of C*-algebras, there are special Cartan subalgebras which are distinguished in the sense that the existence of an arbitrary C*-algebra isomorphism implies the existence of a C*-algebra isomorphism preserving our special Cartan subalgebras. In Section~6, we present such examples by combining rigidity results by Spakula and Willett \cite {SW} and by Whyte \cite{Why} with recent results of the first named author: Under the assumption that the finitely generated groups $\Gamma$ and $\Lambda$ are exact, the isomorphism of their stable uniform Roe algebras implies the isomorphism of their stable groupoids; if moreover $\Gamma$ and $\Lambda$ are non-amenable, this -- and hence in particular isomorphism of uniform Roe algebras -- implies the isomorphism of their groupoids.

%------------------------------------------------

\section{Cartan subalgebras in homogeneous C*-algebras}
\label{sec:hom_gen}

Let $A$ be a $n$-homogeneous C*-algebra over a locally compact Hausdorff space $T$. $A$ is the C*-algebra of continuous sections of a C*-bundle $\cA$ over $T$ with fibre $M_n \defeq M_n(\Cz)$ (see \cite{Dix,DF1,DF2,Fel}). We know that $\cA$ is locally trivial (see for instance \cite[Theorem~3.2]{Fel}). Furthermore, assume that $B$ is a Cartan subalgebra of $A$. Then $B$ is the C*-algebra of continuous sections of a sub-C*-bundle $\cB$ of $\cA$ over the base space $T$. The fibres of $\cB$ are given by the C*-algebra $D_n$, which is the up to conjugacy unique Cartan subalgebra of $M_n$ (see for instance \cite[\S~2.1]{Gre}).

\setlength{\parindent}{0cm} \setlength{\parskip}{0cm}
\blemma
\label{AB:loctriv}
There exists a local trivialization for $(\cA, \cB)$. This means that for every $t \in T$, there exists a neighbourhood $W$ of $t$ and an isomorphism of C*-bundles $\cA \vert_W \overset{\cong}{\lori} W \times M_n$ which identifies $\cB \vert_W$ with $W \times D_n$.
\elemma
\bproof
Fix $t \in T$. By local triviality, there exists an open neighbourhood $V$ of $t$ such that $\cA \vert_V$ is trivial, i.e., we can find an isomorphism $\cA \vert_V \overset{\cong}{\lori} V \times M_n, \, a \ma (s,a(s))$. Using functional calculus, we can find elements $b_1, \dotsc, b_n \in B$ and a neighbourhood $W$ of $t$ with $W \subseteq V$ such that for every $s \in W$, $\gekl{b_1(s), \dotsc, b_n(s)}$ are non-zero pairwise orthogonal projections. Let $\gekl{\zeta_1, \dotsc, \zeta_n}$ be an orthonormal basis of $\Cz^n$ with $b_i(t) \zeta_i = \zeta_i$ for $1 \leq i \leq n$. Then we must have $b_i(t) \zeta_j = \delta_{ij} \zeta_i$ for all $1 \leq i, j \leq n$. By replacing $W$ by a smaller neighbourhood if necessary, we may assume that $b_i(s) \zeta_i \neq 0$ for all $1 \leq i \leq n$ and $s \in W$. Define
$$\eta_i(s) \defeq \frac{b_i(s) \zeta_i}{\norm{b_i(s) \zeta_i}}$$
for all $s \in W$. Since $\gekl{b_1(s), \dotsc, b_n(s)}$ are pairwise orthogonal, $\gekl{\eta_1(s), \dotsc, \eta_n(s)}$ must be an orthonormal basis of $\Cz^n$ for all $s \in W$. Moreover, $b_i(s) \xi = \spkl{\xi,\eta_i(s)}\eta_i(s)$ for all $s \in W$ and $1 \leq i \leq n$. Here $\spkl{\cdot,\cdot}$ denotes the standard scalar product of $\Cz^n$, linear in the first variable. Define the unitary operator $U(s): \: \Cz^n \to \Cz^n, \, (\lambda_i)_i \ma \sum_i \lambda_i \eta_i(s)$. We obtain a C*-bundle isomorphism
$$
  \cA \vert_W \overset{\cong}{\lori} W \times M_n \overset{\cong}{\lori} W \times M_n, \, a \ma (s,a(s)) \ma (s,U(s)^*a(s)U(s)).
$$
Under this isomorphism, $b_i$ is identified with $U(s)^*b_i(s)U(s)$. If $\gekl{e_1, \dotsc, e_n}$ is the standard orthonormal basis of $\Cz^n$, then
$$
  U(s)^*b_i(s)U(s)e_j = U(s)^*b_i(s)\eta_j(s) = \delta_{ij} U(s)^* \eta_i(s) = \delta_{ij} e_i.
$$
Thus, $U(s)^*b_i(s)U(s)$ is the rank one projection corresponding to $e_i$. This shows that $\cB \vert_W$ is identified with $W \times D_n$ under the isomorphism $\cA \vert_W \overset{\cong}{\lori} W \times M_n$ above.
\eproof

In the following, by a local trivialization of $(\cA,\cB)$, we mean a covering of $T$ by open subspaces $W$ together with isomorphisms $\cA \vert_W \overset{\cong}{\lori} W \times M_n$ which identify $\cB \vert_W$ with $W \times D_n$.

\subsection{Reformulation using the language of principal bundles}
\label{ss:PrinBund}

Lemma~\ref{AB:loctriv} allows us to reformulate the existence and uniqueness question for Cartan subalgebras using the language of principal bundles. In the following, let $\Aut(M_n)$ be the automorphism group of the C*-algebra $M_n$, and define
$$
  \Aut(M_n,D_n) \defeq \menge{\varphi \in \Aut(M_n)}{\varphi(D_n) = D_n}.
$$
As explained in \cite[Part~I, \S~8.2]{Ste}, isomorphism classes of locally trivial $M_n$-bundles over $T$ correspond to isomorphism classes of principal $\Aut(M_n)$-bundles over $T$.

\bprop
\label{nhomCartanEx}
An $n$-homogeneous C*-algebra $A$ has a Cartan subalgebra if and only if the principal $\Aut(M_n)$-bundle of the $M_n$-bundle corresponding to $A$ reduces to a principal $\Aut(M_n,D_n)$-bundle.
\eprop
\bproof
Lemma~\ref{AB:loctriv} shows the implication \an{$\Rarr$}: If $\alpha$ is the principal $\Aut(M_n)$-bundle of the $M_n$-bundle $\cA$ corresponding to $A$, then Lemma~\ref{AB:loctriv} shows that there is a principal $\Aut(M_n,D_n)$-bundle $\beta$ over $T$ such that $\alpha$ is isomorphic to the principal $\Aut(M_n)$-bundles over $T$ constructed from $\beta$ and the canonical action $\Aut(M_n,D_n) \curvearrowright \Aut(M_n)$ by left multiplication (see for instance \cite[Chapter~4, \S~5]{Hus} for the construction). This is equivalent to saying that $\alpha$ reduces to a principal $\Aut(M_n,D_n)$-bundle.
\setlength{\parindent}{0.5cm} \setlength{\parskip}{0cm}

Conversely, assume that we can find a principal $\Aut(M_n,D_n)$-bundle $\beta$ over $T$ such that $\alpha$ is isomorphic to the principal $\Aut(M_n)$-bundles over $T$ constructed from $\beta$ and the canonical action $\Aut(M_n,D_n) \curvearrowright \Aut(M_n)$ by left multiplication. Now use $\beta$ to construct an $M_n$-bundle $\cA'$ over $T$. As $\alpha$ is isomorphic to the principal $\Aut(M_n)$-bundles over $T$ constructed from $\beta$ and the canonical action $\Aut(M_n,D_n) \curvearrowright \Aut(M_n)$ by left multiplication, $\cA'$ can be identified with the $M_n$-bundle corresponding to $A$. Moreover, use $\beta$ to construct a $D_n$-bundle $\cB'$ over $T$. $\cB'$ is a sub-C*-bundle of $\cA'$. Therefore, by going over to C*-algebras of continuous sections, we obtain a sub-C*-algebra $B'$ of the C*-algebra $A'$ corresponding to $\cA'$. $B'$ is a Cartan subalgebra of $A'$ because $B'$ and $A'$ correspond to locally trivial C*-bundles by construction, and it is obvious that for an open subspace $W$ of $T$, $C_0(W,D_n)$ is a Cartan subalgebra of $C_0(W,M_n)$.
\eproof
\setlength{\parindent}{0cm} \setlength{\parskip}{0cm}

Given a group $G$ and two principal $G$-bundles $\beta_1$, $\beta_2$ over a base space $T$, we say that $\beta_1$ and $\beta_2$ are isomorphic over $T$ if there is a bundle isomorphism between $\beta_1$ and $\beta_2$ which fixes (every point of) the base space. We say that $\beta_1$ and $\beta_2$ are isomorphic, written $\beta_1 \cong \beta_2$, if $\beta_1$ and $\beta_2$ are isomorphic via a bundle isomorphism which is allowed to induce a non-trivial homeomorphism on $T$. It is the latter notion which corresponds to isomorphism of C*-algebras of continuous sections (see \cite[10.5.5]{Dix}).
\setlength{\parindent}{0cm} \setlength{\parskip}{0.25cm}

Moreover, given two C*-algebras $A_1$ and $A_2$ with Cartan subalgebras $B_1 \subseteq A_1$ and $B_2 \subseteq A_2$, we write $(A_1,B_1) \cong (A_2,B_2)$ if there exists an isomorphism $\varphi: \: A_1 \overset{\cong}{\lori} A_2$ satisfying $\varphi(B_1) = B_2$.

\bprop
\label{nhomCartanUni}
Let $B_1$ and $B_2$ be Cartan subalgebras of an $n$-homogeneous C*-algebra $A$ over $T$. As $B_1$ and $B_2$ are C*-algebras of continuous sections of locally trivial $D_n$-bundles over $T$, they correspond to principal $\Aut(M_n,D_n)$-bundles $\beta_1$ and $\beta_2$. We have $(A,B_1) \cong (A,B_2)$ if and only if $\beta_1 \cong \beta_2$ (as principal $\Aut(M_n,D_n)$-bundles).
\eprop
\bproof
Let $\cA$ be the $M_n$-bundle corresponding to $A$, and $\cB_1$, $\cB_2$ the $D_n$-bundles corresponding to $B_1$, $B_2$. If $(A,B_1) \cong (A,B_2)$, then there exists a bundle isomorphism $\cA \overset{\cong}{\lori} \cA$ (possibly inducing a non-trivial homeomorphism on $T$) which sends $\cB_1$ to $\cB_2$. Looking at local trivializations for $(\cA,\cB_1)$ and $(\cA,\cB_2)$, we see that $\beta_1 \cong \beta_2$ as principal $\Aut(M_n,D_n)$-bundles.
\setlength{\parindent}{0.5cm} \setlength{\parskip}{0cm}

Conversely, if $\beta_1 \cong \beta_2$, then for the $M_n$-bundles $\cA_i$ and $D_n$-bundles $\cB_i$ attached to $\beta_i$ ($i = 1,2$), we obtain an isomorphism $\cA_1 \overset{\cong}{\lori} \cA_2$ sending $\cB_1$ to $\cB_2$. Passing over to the C*-algebras of continuous sections, we get $(A,B_1) \cong (A,B_2)$ because the C*-algebras of continuous sections for $\cA_1$ and $\cA_2$ are both isomorphic to $A$.
\eproof

We summarize our discussion as follows: Fix a base space $T$ (a locally compact Hausdorff space, as above). Consider the canonical map $\cC$ from isomorphism classes of principal $\Aut(M_n,D_n)$-bundles over $T$ to isomorphism classes of principal $\Aut(M_n)$-bundles over $T$, sending the isomorphism class of a principal $\Aut(M_n,D_n)$-bundle $\alpha$ over $T$ to the isomorphism class of the principal $\Aut(M_n)$-bundle constructed from $\alpha$ and the canonical action $\Aut(M_n,D_n) \curvearrowright \Aut(M_n)$ by left multiplication. Here our isomorphisms are allowed to induce non-trivial homeomorphisms on $T$.

By Proposition~\ref{nhomCartanEx}, the $n$-homogeneous C*-algebras over $T$ which have a Cartan subalgebra correspond exactly to the principal $\Aut(M_n)$-bundles which are in the image of $\cC$. Given $[\alpha] \in \img(\cC)$ (here we denote isomorphism classes by $[\cdot]$), let $A$ be the corresponding $n$-homogeneous C*-algebra. Two Cartan subalgebras $B_1$ and $B_2$ in $A$ are called equivalent if $(A,B_1) \cong (A,B_2)$. By Proposition~\ref{nhomCartanUni}, the equivalence classes of Cartan subalgebras in $A$ are in one-to-one correspondence with the pre-image $\cC^{-1}[\alpha]$ of $[\alpha]$.

We can also rephrase this observation using the language of classifying spaces, as long as our base spaces are CW-complexes. Let $G = \Aut(M_n)$ and $H = \Aut(M_n,D_n)$. Clearly, $G \cong U(n) / \Tz$, while $H \cong (\Tz^n \rtimes S_n) / \Tz \cong (\Tz^n / \Tz) \rtimes S_n$. Here $U(n)$ is the (topological) group of unitary $n \times n$-matrices over $\Cz$, and $\Tz$ sits in $U(n)$ as diagonal matrices. This identifies $\Tz$ with the center of $U(n)$. For $H$, $S_n$ is the group of permutations of $\gekl{1, \dotsc, n}$, acting on $\Tz^n$ in the canonical way. This gives rise to the semidirect product $\Tz^n \rtimes S_n$. Again, we embed $\Tz$ into $\Tz^n$, and hence into $\Tz^n \rtimes S_n$ in a canonical way as diagonal matrices. Let $[{\rm Homeo}(T)]$ be the image of the group of homeomorphisms of $T$ in the set of homotopy classes $[T,T]$ of continuous maps from $T$ to $T$. We view $[{\rm Homeo}(T)]$ as a group under composition. Using \cite[Chapter~4, \S~13]{Hus}, the map $\cC$ from above becomes the following map induced by the inclusion $H \into G$:
\bgl
\label{C:HG_T}
  [T,BH] / [{\rm Homeo}(T)] \lori [T,BG] / [{\rm Homeo}(T)]
\egl
where $BH$ and $BG$ are the classifying spaces of $H$ and $G$, and $[{\rm Homeo}(T)$ acts on $[T,BH]$ and $[T,BG]$ by pre-composition. As we have seen above, the discussion of existence and uniqueness of Cartan subalgebras in $n$-homogeneous C*-algebras over $T$ boils down to studying the image and pre-images of the map \eqref{C:HG_T}, in exactly the same way as we explained above.

\subsection{Cartan subalgebras in homogeneous C*-algebras over spheres}

Let us now restrict to the case where our base space $T$ is given by the $k$-sphere
$$S^k = \menge{(x_1, \dotsc, x_{k+1}) \in \Rz^{k+1}}{x_1^2 + \dotso + x_{k+1}^2 = 1}.$$
As above, let $H = \Aut(M_n,D_n)$ and $G = \Aut(M_n)$. By \cite[Chapter~8, Theorem~8.2]{Hus}, isomorphism classes of principal $H$-bundles over $S^k$ are in one-to-one correspondence to $C(S^{k-1},H) / {}_{\sim}$. Here $C(S^{k-1},H)$ is the set of continuous maps $S^{k-1} \to H$, and we define an equivalence relation by saying $c_1 \sim c_2$ if there is $h \in H$ such that $c_1$ and $h c_2 h^{-1}$ are homotopic, for $c_1, c_2 \in C(S^{k-1},H)$. Similarly, isomorphism classes of principal $G$-bundles over $S^k$ are in one-to-one correspondence to $C(S^{k-1},G) / {}_{\sim}$.

As explained in \cite[Chapter~8, Theorem~8.2]{Hus}, these one-to-one correspondences are established by decomposing $S^k$ into upper and lower hemispheres, finding a trivialization of our principal bundle over these hemispheres, and taking the $\sim$-class of the resulting transition function on $S^{k-1}$. Here we embed $S^{k-1}$ into $S^k$ via $(x_1, \dotsc, x_k) \ma (x_1, \dotsc, x_k, 0)$. Thus the map $\cC$ introduced above corresponds to the canonical map $C(S^{k-1},H) / {}_{\sim} \to C(S^{k-1},G) / {}_{\sim}$ induced by the canonical inclusion $H \into G$. Note that all the maps in $C(S^{k-1},H)$ and $C(S^{k-1},G)$ are supposed to send a fixed based point of $S^{k-1}$ to the identity (in $H$ or $G$, respectively). Also note that strictly speaking, \cite{Hus} only discusses isomorphism classes of principal bundles under isomorphisms which induce the identity on the base space. However, every homeomorphism $S^k \overset{\cong}{\lori} S^k$ is homotopic to the identity or to $S^k \to S^k, \, (x_1, \dotsc, x_{k+1}) \ma (x_1, \dotsc, -x_{k+1})$, because the homotopy class is determined by the degree. Both of these maps restrict to the identity on $S^{k-1}$, which we embed into $S^k$ via $S^{k-1} \into S^k, \, (x_1, \dotsc, x_k) \ma (x_1, \dotsc, x_k, 0)$. Therefore, the action of $[{\rm Homeo}(S^k)]$ is trivial, and in our case of spheres, it does not matter whether we allow isomorphisms of principal bundles to induce non-trivial homeomorphisms on the base space or not.

Using this, isomorphism classes of $n$-homogeneous C*-algebras over spheres were studied in \cite{KL}. As explained in \cite{KL}, it turns out that isomorphism classes of principal $G$-bundles over $S^k$ are in one-to-one correspondence to elements in $\pi_{k-1}(G)$. In particular, every $n$-homogeneous C*-algebra over $S^1$ is trivial, i.e., is isomorphic to $C(S^1,M_n)$. Over $S^2$, there are exactly $n$ isomorphism classes of $n$-homogeneous C*-algebras, and an explicit list of such C*-algebras, one from each isomorphism class, is given in \cite{KL}.

In the following discussion, let us use the more precise notation $H_n \defeq \Aut(M_n,D_n)$ and $G_n \defeq \Aut(M_n)$. Our goal is to determine the (number of) isomorphism classes of principal $H_n$-bundles over spheres. We then apply our findings in \S~\ref{ss:PrinBund} to study existence and uniqueness of Cartan subalgebras in $n$-homogeneous C*-algebras over spheres.

Let us start with $k=1$. It is clear that $\abs{C(S^0,H_n) / {}_{\sim}}$ is given by $p(n)$, the number of conjugacy classes of $S_n$, which is the same as the number of partitions of $\gekl{1, \dotsc, n}$. Therefore, we see that the $n$-homogeneous C*-algebra $C(S^1,M_n)$ has exactly $p(n)$ pairwise inequivalent Cartan subalgebras. An explicit construction of these Cartan subalgebras is given as follows: Decompose $S^1$ into $D^1_+ = \menge{(x_1,x_2) \in S^1}{x_2 > -\ve}$ and $D^1_- = \menge{(x_1,x_2) \in S^1}{x_2 < \ve}$, for some fixed $\ve \in (0,\frac{1}{4})$. Choose representatives $\sigma_1, \dotsc, \sigma_{p(n)}$ for the conjugacy classes of $S_n$. For every $1 \leq i \leq p(n)$, we obtain a $D_n$-bundle $\cB_i$ over $S^1$ by gluing together the trivial $D_n$-bundles on $D^1_+$ and $D^1_-$ using the identity on $\menge{(x_1,x_2) \in D^1_+ \cap D^1_-}{x_1 > 0}$ and $\sigma_i$ on $\menge{(x_1,x_2) \in D^1_+ \cap D^1_-}{x_1 < 0}$. These $\cB_i$ give rise to a complete list of pairwise inequivalent Cartan subalgebras $B_1, \dotsc, B_{p(n)}$ of $C(S^1,M_n)$ up to equivalence. For $n=2$, we obtain the trivial bundle and the M{\"o}bius bundle.

Let us consider $k=2$. As we only consider continuous maps in $C(S^1,H_n)$ sending a fixed based point of $S^1$ to the identity of $H_n$, we have $C(S^1,H_n) / {}_{\sim} = C(S^1,\Tz^n / \Tz) / {}_{\sim}$. Here, given $c_1, c_2 \in C(S^1,\Tz^n / \Tz)$, we set $c_1 \sim c_2$ if there exists $\sigma \in S_n$ such that $c_1$ and $\sigma c_2 \sigma^{-1}$ are homotopic. It is then obvious that $C(S^1,\Tz^n / \Tz) / {}_{\sim} \cong \Zz / n \Zz$. Moreover, combining this observation with the discussion in \cite{KL}, it is clear that the map $C(S^1,H_n) / {}_{\sim} \to C(S^1,G_n) / {}_{\sim}$ induced by the canonical inclusion $H_n \into G_n$ is one-to-one. So over $S^2$, every $n$-homogeneous C*-algebra has a Cartan subalgebra, and it is unique up to equivalence.

Finally, we turn to the remaining case $k > 2$. As in the case $k=2$, we have $C(S^{k-1},H_n) / {}_{\sim} = C(S^{k-1},\Tz^n / \Tz) / {}_{\sim}$. Since $\pi_{k-1}(\Tz)$ is trivial for $k > 2$, we conclude that $\abs{C(S^{k-1},\Tz^n / \Tz) / {}_{\sim}} = 1$. This shows that over $S^k$, for $k > 2$, only the trivial $n$-homogeneous C*-algebra $C(S^k,M_n)$ has a Cartan subalgebra. Moreover, up to equivalence, $C(S^k,M_n)$ has a unique Cartan subalgebra given by $C(S^k,D_n)$. This gives a conceptual explanation for Gregson's and Natsume's examples of homogeneous C*-algebras over spheres which do not admit Cartan subalgebras (see \cite[\S~2.2.5]{Gre} and the appendix of \cite{Kum}).

Exactly the same methods allow us to study Cartan subalgebras of $\aleph_0$-homogeneous C*-algebras over spheres, at least those where the Cartan pair admits a local trivialization. First of all, fix a countably infinite dimensional Hilbert space $\cH$ and let $\cK \defeq \cK(\cH)$ be the C*-algebra of compact operators on $\cH$. $\aleph_0$-homogeneous C*-algebras over spheres correspond to locally trivial $\cK$-bundles, and two such $\aleph_0$-homogeneous C*-algebras are isomorphic if their Dixmier-Douady invariants coincide (see for instance \cite[Theorem~IV.1.7.15]{Bla}). We know that $H^3(S^k,\Zz) \cong \gekl{0}$ if $k \neq 3$ and $H^3(S^3,\Zz) \cong \Zz$. Hence over $S^k$, for $k \neq 3$, every $\aleph_0$-homogeneous C*-algebra is trivial, i.e., isomorphic to $C(S^k,\cK)$, while there are countably infinitely many pairwise non-isomorphic $\aleph_0$-homogeneous C*-algebras over $S^3$. Moreover, by \cite{R1983} (see also \cite[\S~6.4]{R08}), every $\aleph_0$-homogeneous C*-algebra over $S^k$ has a Cartan subalgebra. Since  $\aleph_0$-homogeneous C*-algebras $A$ over a finite dimensional second countable Hausdorff locally compact space $T$ are classified up to isomorphism by their Dixmier-Douady invariant $\delta(A)\in H^3(T,\Zz)$, it suffices to construct a groupoid model for such an algebra. We realize the class $\delta(A)$ in $H^3(T,\Zz)$ as a \v Cech cocycle $(\sigma_{ijk})$ on an open cover $(U_i)_{i\in I}$ of $T$; thus $\sigma_{ijk}$ is a continuous $\Tz$-valued function defined on $U_i\cap U_j\cap U_k$. In order to get infinite fibers, we replace the index set $I$ by $\tilde I:I\times\Nz$ and write $U_{(i,n)}:=U_i$. We let $X:=\{(\tilde i,t)\in \tilde I\times T: t\in U_{\tilde i}\}$ be the disjoint union of the $U_{\tilde i}$s and $\pi: X\to T$ be the second projection. It is a local homeomorphism with infinite fibers. The associated \'etale equivalence relation $R=X\times_T X$ carries the twist defined by the 2-cocycle $\sigma\in Z^2(R,\Tz)$ defined by 
$\sigma((\tilde i,t),(\tilde j,t),(\tilde k,t))=\sigma_{ijk}(t)$. 
The twisted groupoid C*-algebra $C^*(R,\sigma)$ has $C(X)$ as a Cartan subalgebra. Since $\pi: X\to T$ is usually not a covering map, the corresponding Cartan pair usually does not admit local trivializations. Also, another choice of open cover may give another Cartan pair. However, the twisted groupoids obtained from this construction are all equivalent.
\setlength{\parindent}{0.5cm} \setlength{\parskip}{0cm}

Note that the C*-algebras attached to groupoids twisted by a 2-cocycle are a special case of C*-algebras of twisted groupoids. \cite[Chapter~I, Proposition~1.14]{Ren} explains how to construct the twist $\Sigma=E_\sigma$ from a 2-cocycle $\sigma$. In that case, $\Sigma$ is a trivial principal $\Tz$-bundle. In fact, a twist $\Sigma$ is of the form $E_\sigma$ if and only if it is a trivial principal $\Tz$-bundle. As an example of a twist which is not given by a 2-cocycle, just consider a non-trivial principal $\Tz$-bundle over a space $G=X$.
\setlength{\parindent}{0cm} \setlength{\parskip}{0.25cm}

Let us study existence and uniqueness of Cartan subalgebras of $\aleph_0$-homogeneous C*-algebras where the corresponding Cartan pairs admit local trivializations. Let $\cD$ be a fixed Cartan subalgebra of $\cK$. It does not matter which one we choose as every two such Cartan subalgebras are equivalent. Given an $\aleph_0$-homogeneous C*-algebra $A$ over $T$ and a Cartan subalgebra $B$, we let $\cA$ and $\cB$ be the corresponding $\cK$- and $\cD$-bundles over $T$. Recall that a local trivialization of $(\cA,\cB)$ is a covering of $T$ by open subsets $W$ together with isomorphisms $\cA \vert_W \overset{\cong}{\lori} W \times \cK$ sending $\cB \vert_W$ to $W \times \cD$. Let $G = \Aut(\cK) \cong U(\cH) / \Tz$ and $H = \Aut(\cK,\cD) = (\Tz^{\infty} / \Tz) \rtimes S_{\infty}$. Here $S_{\infty}$ is the group of permutations of $\gekl{1, 2, \dotsc}$. The same argument as in the finite homogeneous case shows that over $T = S^k$, $\aleph_0$-homogeneous C*-algebras are in one-to-one correspondence with $C(S^{k-1},G) / {}_{\sim}$. Moreover, an $\aleph_0$-homogeneous C*-algebra has a Cartan subalgebra such that the corresponding Cartan pair has a local trivialization if and only if it corresponds to an element in $C(S^{k-1},G) / {}_{\sim}$ which lies in the image of the canonical map $C(S^{k-1},H) / {}_{\sim} \to C(S^{k-1},G) / {}_{\sim}$. Given such an element in $C(S^{k-1},G) / {}_{\sim}$, equivalence classes of Cartan subalgebras in the corresponding $\aleph_0$-homogeneous C*-algebra are in one-to-one correspondence with the pre-image of that element under the canonical map $C(S^{k-1},H) / {}_{\sim} \to C(S^{k-1},G) / {}_{\sim}$.

For $k=1$, it is easy to see that $C(S^{k-1},H) / {}_{\sim}$ is in one-to-one correspondence with conjugacy classes in $S_{\infty}$, and there are uncountably many of these. For $k \geq 2$, we have $\abs{C(S^{k-1},H) / {}_{\sim}} = 1$.

This shows that for $k=1$, the up to isomorphism unique $\aleph_0$-homogeneous C*-algebra $C(S^1,\cK)$ over $S^1$ has uncountably many Cartan subalgebras whose Cartan pairs admit local trivializations. For $k=2$ or $k \geq 4$, we see that the up to isomorphism unique $\aleph_0$-homogeneous C*-algebra $C(S^k,\cK)$ over $S^k$ has $C(S^k,\cD)$ as its up to equivalence unique Cartan subalgebra whose Cartan pair admit a local trivialization. For $k=3$, the only $\aleph_0$-homogeneous C*-algebra over $S^3$ which admits a Cartan subalgebra whose Cartan pair has a local trivialization is the trivial one, i.e., isomorphic to $C(S^3,\cK)$. And up to equivalence, the Cartan subalgebra of $C(S^3,\cK)$ whose Cartan pair admit a local trivialization is given by $C(S^3,\cD)$. All the other non-trivial $\aleph_0$-homogeneous C*-algebra over $S^3$ also have Cartan subalgebras, as we mentioned above. Therefore, these Cartan subalgebras do not have the property that the corresponding Cartan pairs admit local trivializations.

The case $k=3$ shows that not every Cartan subalgebra of an $\aleph_0$-homogeneous C*-algebra is such that the corresponding Cartan pair admits a local trivialization. It is also interesting to point out that for $k=4$, by \cite{KL}, there exist non-trivial $n$-homogeneous C*-algebras over $S^4$, for every $n \geq 2$. By our previous discussion, these non-trivial $n$-homogeneous C*-algebras cannot have Cartan subalgebras. However, after stabilizing, i.e., tensoring with $\cK$, we obtain the trivial $\aleph_0$-homogeneous C*-algebra over $S^4$, because its Dixmier-Douady invariant has to vanish. And this trivial $\aleph_0$-homogeneous C*-algebra obviously has a Cartan subalgebra. Thus we have just obtained an example of a C*-algebra which has no Cartan subalgebra, but whose stabilization admits a Cartan subalgebra.

\section{Cartan subalgebras in group C*-algebras of certain connected Lie groups}

The structure of reduced and full group C*-algebras of Lie groups have been studied by many authors (see for instance \cite{Fel,Mil,KM,Val,PP}). Our aim is to show that many of these group C*-algebras have a Cartan subalgebra.

We start with the full group C*-algebra $C^*(SL(2,\Cz))$ of $SL(2,\Cz)$. Its structure has been determined in \cite{Fel}. Recall that by \cite[\S~5]{Fel}, the C*-algebra $C^*(SL(2,\Cz))$ can be described as follows: Let $Z_1 \defeq \menge{(m,\rho)}{m \in \Zz_{\geq 0}, \, \rho \in \Rz, \, \rho \geq 0 \ {\rm if} \ m = 0}$, and endow $Z_1$ with the subspace topology of $\Zz \times \Rz$. Let $Z_2 = [0,1]$ with the usual topology. Define $Z$ to be the topological space obtained from the disjoint union of $Z_1$ and $Z_2$ by identifying $(0,0) \in Z_1$ with $0 \in Z_2$. Moreover, let $H = \ell^2(\Zz_{\geq 0})$ with canonical orthonormal basis $e_0, e_1, e_2, \dotsc$, and let $K$ be the closed subspace of $H$ generated by $\menge{e_n}{n \geq 1}$. Let $M: \: K \to H$ be the unitary given by $M(e_n) = e_{n-1}$ for all $n \geq 1$. We write $\cK(H)$ for the C*-algebra of compact operators on $H$. Let $A$ be the C*-algebra of all $a \in C_0(Z,\cK(H))$ satisfying the condition that
$$a(1) \in \cK(\Cz e_0) \oplus M^{-1} a(2,0) M \subseteq \cK(\Cz e_0) \oplus \cK(K) \subseteq \cK(H).$$
Then \cite[Theorem~5.4]{Fel} tells us that $C^*(SL(2,\Cz)) \cong A$.

Let $\cD$ be the C*-algebra of diagonal compact operators with respect to $\menge{e_n}{n \in \Zz_{\geq 0}}$, i.e., $\cD = \menge{T \in \cK(H)}{T e_n \in \Cz e_n \ {\rm for} \ {\rm all} \ n \in \Zz_{\geq 0}}$.
\bprop
\label{PROP:SL2C}
$B \defeq A \cap C_0(Z,\cD)$ is a Cartan subalgebra in $A$. In particular, $C^*(SL(2,\Cz))$ admits a Cartan subalgebra.
\eprop
\setlength{\parindent}{0cm} \setlength{\parskip}{0cm}

\bproof
Since $C_0(Z,\cD)$ is obviously a Cartan subalgebra of $C_0(Z,\cK(H))$, it is clear that $B$ contains an approximate unit for $A$, that $B$ is maximal abelian in $A$, and that there exists a faithful conditional expectation $A \onto B$. It remains to show that the normalizer
$$
  N_A(B) = \menge{n \in A}{nBn^* \subseteq B \ {\rm and} \ n^*Bn \subseteq B}
$$
generates $A$ as a C*-algebra.

For $g \in C_0(\Rz)$, $f_0 \in C_0(\Rz)$ with $f_0(-1) = 1$ and $f_2 \in C_0(\Rz)$ with $f_2(2) = 1$, consider the function $n: \: Z \to \cK(H)$ given by
\bglnoz
  n(m,\rho) &=& \delta_{m,2} f_2(\rho) e_{i,j} \ {\rm for} \ (m,\rho) \in Z_1, \, m > 0 \\
  n(0,\rho) &=& (g(\rho), f_0(\rho) e_{i+1,j+1}) \in \cK(\Cz e_0) \oplus \cK(K) \ {\rm for} \ (0,\rho) \in Z_1 \\
  n(t) &=& (g(-t), f_0(-t) e_{i+1,j+1}) \in \cK(\Cz e_0) \oplus \cK(K) \ {\rm for} \ t \in Z_2.
\eglnoz
Here $\gekl{e_{i,j}}_{i,j}$ are the canonical matrix units with respect to the orthonormal basis $\gekl{e_n}_n$.
\setlength{\parindent}{0.5cm} \setlength{\parskip}{0cm}

It is clear that $n \in A$ and $n \in N_A(B)$. Moreover, given $a \in A$ arbitrary, we can find elements $a_i \in A$ in the linear span of elements $n$ as above such that $\lim_{i \to \infty} (a-a_i)(2,0) = 0$ and $\lim_{i \to \infty} (a-a_i)(1) = 0$. Therefore, it suffices to show that every $a \in A$ with $a(2,0) = a(1) = 0$ lies in the C*-algebra generated by $N_A(B)$. The collection of these elements obviously forms the ideal $I \defeq C_0(Z \setminus \gekl{(2,0),1}, \cK(H))$. We have $I \cap B = I \cap C_0(Z,\cD) = C_0(Z \setminus \gekl{(2,0),1}, \cD)$. It is clear that $I \cap B$ is a Cartan subalgebra of $I$. In particular, $I$ is generated as a C*-algebra by $N_I(I \cap B)$. Hence all we have to show is $N_I(I \cap B) \subseteq N_A(B)$. Choose $h \in N_I(I \cap B)$ and $b \in B$. Moreover, let $h_{\lambda} \in I \cap B$ be an approximate unit for $I$. Then $nbn^* = \lim_{\lambda} n h_{\lambda} b h_{\lambda} n^* \in I \cap B$ as $h_{\lambda} b h_{\lambda} \in I \cap B$. Thus, indeed, $N_I(I \cap B) \subseteq N_A(B)$.
\eproof

To cover more group C*-algebras of Lie groups, we prove the following general result: Let $M$ be a locally compact Hausdorff space. Let $H = \ell^2 \Zz$ and $\menge{e_i}{i \in \Zz}$ be the canonical orthonormal basis of $H$. Let $\cK(H)$ be the C*-algebra of compact operators on $H$ and $\cD$ the C*-algebra of compact operators on $H$ which are diagonal with respect to $\menge{e_i}{i \in \Zz}$. Assume that for every $p \in M$, we are given a decomposition $H = H_{1,p} \oplus \dotso \oplus H_{k(p),p}$. Moreover, suppose that there exists an open covering $\cV$ of $M$ such that for every $V \in \cV$, there is a continuous map $V \to U(H), \, p \ma U^V_p$ and 
\begin{align}
\label{c1}
  & {\rm there} \ {\rm is} \ {\rm a} \ {\rm partition} \ \Zz = P_{1,p} \sqcup \dotso \sqcup P_{k(p),p} \ {\rm for} \ {\rm every} \ p \in M \ {\rm with} \\
  & H_{i,p} = \spkl{\menge{U^V_p(e_j)}{j \in P_{i,p}}} \nonumber
\end{align}
for all $1 \leq i \leq k(p)$, and
%\label{c1}
 % {\rm there} \ {\rm is} \ {\rm a} \ {\rm partition} \ \Zz = P_1 \sqcup \dotso \sqcup P_{k(p)} \ {\rm for} \ {\rm every} \ p \in M \ {\rm with} \ H_{i,p} = \spkl{\menge{U^V_p(e_j)}{j \in P_i}}
%\egl
%for all $1 \leq i \leq k(p)$, and
%\bgl
\bgl
\label{c2}
  \menge{U^V_p(e_i)}{i \in \Zz} = \menge{U^W_p(e_i)}{i \in \Zz}
\egl
for all $V, W \in \cV$ and $p \in V \cap W$.
The latter condition allows us to define, for every $p \in M$, $\cD_p \defeq U^V_p \cD (U^V_p)^*$, where $V \in \cV$ satisfies $p \in V$ (by \eqref{c2}, $\cD_p$ does not depend on the choice of $V$). Define
$$
  A \defeq \menge{a \in C_0(M,\cK(H))}{a(p) \in \cK(H_{1,p}) \oplus \dotso \oplus \cK(H_{k(p),p}) \ {\rm for} \ {\rm all} \ p \in M},
$$
$$
  B \defeq \menge{b \in A}{b(p) \in \cD_p \ {\rm for} \ {\rm all} \ p \in M}.
$$
\bprop
\label{Prop:genLie}
$B$ is a Cartan subalgebra of $A$.
\eprop
\bproof
Our conditions allow us to work locally. And working locally, all we have to show is that given an open subset $V$ of $M$, and a sub-C*-algebra $A_V$ of $C_0(V,\cK(H))$ containing $B_V \defeq C_0(V,\cD)$, $B_V$ is a Cartan subalgebra of $A_V$.
\setlength{\parindent}{0.5cm} \setlength{\parskip}{0cm}

It is clear that $B_V$ has an approximate unit for $A_V$, that $B_V$ is maximal abelian in $A_V$, and that there is a faithful conditional expectation $A_V \onto B_V$. All we have to show is that the normalizer $N_{A_V}(B_V) = \menge{n \in A_V}{nB_Vn^* \subseteq B_V \ {\rm and} \ n^*B_Vn \subseteq B_V}$ generates $A_V$ as a C*-algebra. Given an arbitrary element $a \in A_V$, let $e_{i,j}$ be the canonical matrix units attached to the orthonormal basis $\menge{e_i}{i \in \Zz}$, and let $d_i \in C_b(V,\cD)$ be given by $d_i(p) = e_{i,i}$ for all $p \in V$. Obviously, $a = \lim_{n \to \infty} \sum_{i,j=1}^n d_i a d_j$, and
%clear that for all $i,j \in \Zz$, $e_i a e_j$ lies in $N_{A_V}(B_V)$.
$d_i a d_j$ clearly lies in $N_{A_V}(B_V)$ for all $i,j \in \Zz$.
\eproof

Using this, let us describe a Cartan subalgebra in $C^*(SL(2,\Rz))$. The full C*-algebra of $SL(2,\Rz)$ has been described in \cite{Mil}. We recall the description: Define
\begin{align*}
  K \defeq & \ \ \ \menge{(x,0) \in \Cz^2}{x \in [0,\infty) \cdot i} \\
  & \cup \menge{(y,1) \in \Cz^2}{y \in [0,\infty) \cdot i} \\
  & \cup \menge{(z,0) \in \Cz^2}{z \in [0,1]}
\end{align*}
equipped with the subspace topology of $\Cz^2$, and
$$
  M \defeq K \sqcup \Zz_{\geq 2} \sqcup \Zz_{\geq 2},
$$
where $\Zz_{\geq 2}$ is endowed with the discrete topology. As above, let $H = \ell^2 \Zz$. We use the same notation as above. For all $p \in M$, $p \neq (0,1)$, $p \neq (1,0)$, let $k(p) = 1$ and $H_{1,p} = H$. For $p=(0,1)$, let $k(p) = 2$ and $H_{1,p} = \ell^2 (\Zz_{<0}) \subseteq H$ and $H_{2,p} = \ell^2 (\Zz_{\geq 0}) \subseteq H$. For $p = (1,0)$, set $k(p) \defeq 3$ and $H_{1,p} = \ell^2 (\Zz_{<0})$, $H_{2,p} = \Cz e_0$ and $H_{3,p} = \ell^2 (\Zz_{>0})$. Then $A \defeq \menge{a \in C_0(M,\cK(H))}{a(p) \in \cK(H_{1,p}) \oplus \dotso \oplus \cK(H_{k(p),p}) \ {\rm for} \ {\rm all} \ p \in M}$ is isomorphic to $C^*(SL(2,\Rz))$ by \cite{Mil}.
\setlength{\parindent}{0cm} \setlength{\parskip}{0.25cm}

To show that we can apply Proposition~\ref{Prop:genLie}, let $\cV = \gekl{M}$ and $U^M_p \defeq \id_H$. Obviously, the conditions for Proposition~\ref{Prop:genLie} are satisfied. Let $\cD$ be the C*-algebra of compact operators on $H$ which are diagonal with respect to $\menge{e_i}{i \in \Zz}$. We obtain
\bcor
\label{COR:Sl2R}
$C_0(M,\cD)$ is a Cartan subalgebra of $A$. In particular, $C^*(SL(2,\Rz))$ has a Cartan subalgebra.
\ecor
\setlength{\parindent}{0cm} \setlength{\parskip}{0cm}

Let us describe a Cartan subalgebra of the full group C*-algebra $C^*(G)$ of the universal covering group $G$ of $SL(2,\Rz)$. Recall the description of $C^*(G)$ in \cite{KM}: Let $M = \menge{z \in \Cz}{0 < \abs{z} \leq 1} \cup [1,2]$ equipped with the subspace topology of $\Cz$. Let $H = \ell^2 \Zz$ and $U$ the unilateral shift. Choose a self-adjoint operator $T \in \cL(H)$ such that $U = e^{2 \pi i T}$. Define $U: \: [0,1] \to U(H), \, t \ma e^{2 \pi i t T}$. Clearly, $U(0) = \id_H$ and $U(1) = U$. For every $t \in [0,1]$, set $H_t \defeq U(t)^*(\ell^2(\Zz_{\geq 0}))$. Set $k(p) = 1$ if $p \in M$, $0 < \abs{p} < 1$, and $H_{1,p} = H$. Set $k(p) = 2$ if $p \in M$, $\abs{p} = 1$, $p \neq 1$, and $H_{1,p} = H_t^{\perp}$, $H_{2,p} = H_t$, where $p = e^{2 \pi i t}$ for a uniquely determined $t \in (0,1)$. Set $k(p) = 3$ if $p \in [1,2] \subseteq M$, and $H_{1,p} \defeq \ell^2(\Zz_{<0})$, $H_{2,p} \defeq \Cz e_0$, $H_{3,p} \defeq \ell^2(\Zz_{>0})$. Let
$$
  A \defeq \menge{a \in C_0(M,\cK(H))}{a(p) \in \cK(H_{1,p}) \oplus \dotso \oplus cK(H_{k(p),p})},
$$
$$
  A_1 \defeq \menge{a \in A}{a(p) e_0 = a(1) e_0 \ {\rm for} \ {\rm all} \ p \in [1,2] \ {\rm and} \ a(2) \in \Cz e_{0,0}}.
$$
Then $A_1 \cong C^*(G)$ by \cite{KM}.
\setlength{\parindent}{0cm} \setlength{\parskip}{0.25cm}

Let
$$
  V \defeq \menge{r e^{2 \pi i t}}{r \in (0,1], \, t \in (-\frac{\pi}{4},\frac{\pi}{4})} \cup [1,2]
$$
and
$$
  W \defeq \menge{r e^{2 \pi i t}}{r \in (0,1], \, t \in (0,1)}.
$$
It is clear that $M = V \cup W$. Define $U^V_p \defeq e^{-2 \pi i t T}$ if $p = r e^{2 \pi i t}$ ($t \in (-\frac{\pi}{4},\frac{\pi}{4})$), $U^V_p \defeq \id_H$ if $p \in [1,2]$ and $U^W_p \defeq U(t)^*$ if $p = r e^{2 \pi i t}$ ($t \in (0,1)$). In this way, we obviously obtain continuous maps $V \to U(H), \, p \ma U^V_p$ and $W \to U(H), \, p \ma U^W_p$.

Let us check condition \eqref{c1}: For $p \in M$ with $0 < \abs{p} < 1$, condition \eqref{c1} obviously holds. For $p \in M$ with $p = e^{2 \pi i t}$, condition \eqref{c1} is clearly satisfied for $U^W$. For $U^V$ and $p = e^{2 \pi i t}$ with $t \in [0,\frac{\pi}{4})$, it is also obvious that \eqref{c1} holds. For $U^V$ and $p = e^{2 \pi i t}$ with $t \in (-\frac{\pi}{4},0)$, we have
\bglnoz
  H_{2,p} &=& U(1+t)^* \ell^2(\Zz_{\geq 0}) 
  = e^{-2 \pi i (1+t) T} \ell^2(\Zz_{\geq 0})
  = e^{-2 \pi i t T} e^{-2 \pi i T} \ell^2(\Zz_{\geq 0}) \\
  &=& U^V_p U^* \ell^2(\Zz_{\geq 0})
  = U^V_p \ell^2(\Zz_{\geq -1})
\eglnoz
and $H_{1,p} = H_{2,p}^{\perp} = U^V_p \ell^2(\Zz_{\leq -2})$. Hence \eqref{c1} holds. For $p \in [1,2] \subseteq M$, it is again obvious that condition \eqref{c1} is true.

Let us check condition \eqref{c2}. We have
$$
  V \cap W = \menge{r e^{2 \pi i t} \in \Cz}{r \in (0,1], \, t \in (-\frac{\pi}{4},0) \cup (0,\frac{\pi}{4})}.
$$
Clearly, \eqref{c2} holds for $p \in V \cap W$ of the form $p = r e^{2 \pi i t}$ with $ r \in (0,1]$ and $t \in (0,\frac{\pi}{4})$. For $p = r e^{2 \pi i t}$ with $ r \in (0,1]$ and $t \in (-\frac{\pi}{4},0)$, we have
\bglnoz
  && \menge{U^V_p(e_i}{i \in \Zz} = \menge{e^{- 2 \pi i t T}(e_i)}{i \in \Zz} = \menge{e^{- 2 \pi i t T}(e_{i-1})}{i \in \Zz} \\
  &=& \menge{e^{- 2 \pi i t T} e^{- 2 \pi i T} (e_i)}{i \in \Zz}
  = \menge{e^{- 2 \pi i (t+1) T}(e_i)}{i \in \Zz} = \menge{U^W_p(e_i)}{i \in \Zz}.
\eglnoz
Hence condition \eqref{c2} holds as well.

Therefore, Proposition~\ref{Prop:genLie} applies and tells us that if we set $\cD_p \defeq U^V_p \cD (U^V_p)^*$ for $p \in V$ and $\cD_p \defeq U^W_p \cD (U^W_p)^*$ for $p \in W$, then $B \defeq \menge{b \in A}{b(p) \in \cD_p \ {\rm for} \ {\rm all} \ p \in M}$ is a Cartan subalgebra of $A$. Using this, it is now straightforward to check, using similar methods as in the proof of Proposition~\ref{PROP:SL2C}, that $B \cap A_1$ is a Cartan subalgebra of $A_1$. Thus, we obtain
\bcor
The full C*-algebra of the universal covering group of $SL(2,\Rz)$ has a Cartan subalgebra.
\ecor
\setlength{\parindent}{0cm} \setlength{\parskip}{0cm}

Using the description of full group C*-algebras for covering groups with finite center of $SO_0(2,1)$ in \cite{KM}, the same explanation as for Corollary~\ref{COR:Sl2R} gives
\bcor
The full group C*-algebra of a covering group with finite center of $SO_0(2,1)$ has a Cartan subalgebra.
\ecor

\bremark
Looking at the conjecture at the end of \cite{KM}, it seems that our method can be applied to show that for many connected semisimple Lie groups, the full group C*-algebras have Cartan subalgebras.
\eremark

We now turn to reduced group C*-algebras of certain Lie groups. We use the descriptions of these C*-algebras in \cite{BM,Val,PP}.
\bcor
For every connected semi-simple Lie group with real rank one and finite center, the reduced group C*-algebra has a Cartan subalgebra.
\ecor
\bproof
By the description of these reduced group C*-algebras in \cite{BM,Val}, we see that they are of the form as in Proposition~\ref{Prop:genLie}, but with $k(p)$ equal to $1$ or $2$ for all $p \in M$, and if $k(p) = 2$, we can always choose $H_{1,p} = \ell^2(\Zz_{\leq -1})$ and $H_{2,p} = \ell^2(\Zz_{\geq 0})$. Thus, Proposition~\ref{Prop:genLie} applies because we can simply choose the trivial cover $\cV = \gekl{M}$ and the constant map $U_M \equiv \id_H$.
\eproof

\bcor[Corollary to Proposition~4.1 in \cite{PP}]
For every connected, complex semi-simple Lie group, the reduced group C*-algebra has a Cartan subalgebra.
\ecor
\bproof
Because of \cite[Proposition~4.1]{PP}, these reduced group C*-algebras are trivial $\aleph_0$-homogeneous C*-algebras. Hence they obviously admit Cartan subalgebras.
\eproof

\bquestion
Is there a conceptual explanation why the group C*-algebras of (many) connected Lie groups have Cartan subalgebras?
\equestion

\section{Non-existence of C*-diagonals and Cartan subalgebras in certain reduced group C*-algebras}

Using non-existence results for Cartan subalgebras in certain group vN-algebras \cite{Voi,Hay,Oza,OP10a,OP10b,PV14a,PV14b,BHV}, we show that certain reduced group C*-algebras have no diagonals, or even stronger, no Cartan subalgebras.
\setlength{\parindent}{0cm} \setlength{\parskip}{0.25cm}

Let $G$ be an {\'e}tale second countable locally compact Hausdorff groupoid, and suppose that we are given a twist $\Tz \times G^{(0)} \overset{i}{\rightarrowtail} \Sigma \overset{j}{\twoheadrightarrow} G$. As explained in \cite{R08}, every Cartan pair is of the form $(C^*_r(G,\Sigma), C_0(G^{(0)}))$, where $G$ is in addition assumed to be topologically principal. For the notion of twists and the construction of $C^*_r(G,\Sigma)$, we refer to \cite{R08} and the references therein.

\subsection{Traces on twisted groupoid C*-algebras}

Let $X \defeq G^{(0)}$, and for simplicity, let us assume that $X$ is compact. Let $S \subseteq G$ be an open bisection, i.e., the range and source maps restrict to homeomorphisms $r \vert_S: \: S \to r(S)$ and $s \vert_S: \: S \to s(S)$. Let $\alpha_S$ be the homeomorphism $s(S) \to r(S), \, x \ma r((s \vert_S)^{-1}(x))$. It is easily checked that a Borel measure $\mu$ on $X$ is invariant in the sense of \cite[Definition~I.3.12]{Ren} if and only if for every open bisection $S$ in $G$ and every $f \in C_0(r(S))$, we have $\int_{r(S)} f \, {\rm d} \mu = \int_{s(S)} f \circ \alpha_S \, {\rm d} \mu$. In other words, we have $\mu \vert_{r(S)} = \alpha_S (\mu \vert_{s(S)})$. 
%A Borel measure $\mu$ on $X$ is called $\cS(G)$-invariant if for every open bisection $S$ in $G$ and every $f \in C_0(r(S))$, we have $\int_{r(S)} f \, {\rm d} \mu = \int_{s(S)} f \circ \alpha_S \, {\rm d} \mu$. In other words, we require that $\mu \vert_{r(S)} = \alpha_S (\mu \vert_{s(S)})$. 
Now let $\tau$ be a trace on $C^*_r(G,\Sigma)$. $\tau \vert_{C(X)}$ corresponds to a probability measure $\mu$ on $X$. The following lemma is certainly well-known, but we include this result for completeness and because we were not able to locate a reference.
\blemma
\label{muSGinv}
$\mu$ is invariant.
%$\mu$ is $\cS(G)$-invariant.
\elemma
\setlength{\parindent}{0cm} \setlength{\parskip}{0cm}

\bproof
Let $f \in C_c(s(S))$ be positive. Define $g \in C_c(j^{-1}(S))$ by $g(\zeta) = f^{\frac{1}{2}}(r(j(\zeta))$. Then $g * g^*, g^* * g \in C(X)$. For $x \in s(S)$, let $\gamma \in S$ be determined by $s(\gamma) = x$. We have $g * g^* (x) = \rukl{f^{\frac{1}{2}}(r(\gamma^{-1}))}^2 = f(s(\gamma)) = f(x)$. At the same time, since $\alpha_S(x) = r(\gamma)$, we get $g^* * g (x) = \rukl{f^{\frac{1}{2}}(r(\gamma))}^2 = f(r(\gamma)) = f(\alpha_S(x))$. Thus, as $\tau$ is a trace, $\int_{r(S)} f \, {\rm d} \mu = \tau(f) = \tau(g * g^*) = \tau(g^* * g) = \tau(f \circ \alpha_S) = \int_{s(S)} f \circ \alpha_S \, {\rm d} \mu$.
\eproof

Now let us assume that, in addition to the assumptions above, $G$ is also topologically principal. Let $E: \: C^*_r(G,\Sigma) \onto C(X)$ be the canonical faithful conditional expectation. Moreover, let $\mu$ be an invariant probability measure on $X$, viewed as a state on $C(X)$.
%Moreover, let $\mu$ be an $\cS(G)$-invariant probability measure on $X$, viewed as a state on $C(X)$.
\setlength{\parindent}{0cm} \setlength{\parskip}{0.25cm}

%\blue{
%In the next two lemmas: How are cocycles $\sigma$ and twists $\Sigma$ related? Are these notions equivalent, or is one (potentially) more general than the other? (Maybe twists are more general?)
%}
%
%\red{The construction of the twist $\Sigma=E_\sigma$ from the 2-cocycle $\sigma$ is given in [47, I, Proposition 1.14]. In that case, $\Sigma$ is a trivial principal $S^1$-bundle. In fact, a twist $\Sigma$ is of the form $E_\sigma$ if and only if it is a trivial principal $S^1$-bundle. As an example of a twist which is not given by a 2-cocycle, just consider a non trivial principal $S^1$-bundle over a space $G=X$.}

\blemma (see also \cite[Proposition~II.5.4]{Ren})
\label{muE:trace}
$\mu \circ E$ is a trace on $C^*_r(G,\Sigma)$.
\elemma
%\blemma
%\label{muE:trace}
%$\mu \circ E$ is a trace on $C^*_r(G,\Sigma)$.
%\elemma
\setlength{\parindent}{0cm} \setlength{\parskip}{0cm}

\bproof
As in \cite{R08}, we denote by $C_c(G,\Sigma)$ the set of compactly supported, continuous functions $f: \: \Sigma \to \Cz$ such that $f(z \zeta) = f(\zeta) \bar{z}$.
Given $f \in C_c(G,\Sigma)$, we set $\supp(f) \defeq \menge{\gamma \in G}{f(\zeta) \neq 0 \ {\rm for} \ {\rm all} \ \zeta \in \Sigma \ {\rm with} \ j(\zeta) = \gamma}$. Now, to prove the lemma, it suffices to show that $(\mu \circ E) (f * g) = (\mu \circ E) (g * f)$ for all $f, g \in C_c(G,\Sigma)$ such that $\supp(f) \subseteq S$ and $\supp(g) \subseteq T$ for some open bisections $S$ and $T$. For $x \in X$, we have $(f * g)(x) = f(\zeta^{-1}) g(\zeta)$ if $x \in s(S^{-1} \cap T)$ and $\zeta \in \Sigma$ satisfies $s(j(\zeta)) = x$, $j(\zeta) \in S^{-1} \cap T$, and $(f * g)(x) = 0$ if $x \notin s(S^{-1} \cap T)$. Similarly, for $y \in X$, we have $(g * f)(y) = g(\eta^{-1}) f(\eta)$ if $y \in s(S \cap T^{-1})$ and $\eta \in \Sigma$ satisfies $s(j(\eta)) = y$, $j(\zeta) \in S \cap T^{-1}$, and $(g * f)(y) = 0$ if $y \notin s(S \cap T^{-1})$. Thus $(f * g) \vert_X \in C_0(s(S^{-1} \cap T))$ and $(g * f) \vert_X \in C_0(s(S \cap T^{-1})) = C_0(r(S^{-1} \cap T))$. Moreover, for $x \in s(S^{-1} \cap T)$, choose $\zeta \in \Sigma$ with $s(j(\zeta)) = x$, $j(\zeta) \in S^{-1} \cap T$. Then $\alpha_{S^{-1} \cap T}(x) = r(j(\zeta)) = s(j(\zeta^{-1}))$ and $j(\zeta^{-1}) \in S \cap T^{-1}$. Hence $(f * g)(x) = g(\zeta) f(\zeta^{-1}) = (g * f)(s(j(\zeta^{-1})) = (g * f)(\alpha_{S^{-1} \cap T}(x))$. We conclude that
$(\mu \circ E) (f * g) = \mu((f * g) \vert_X) = \mu((g * f) \vert_X \circ \alpha_{S^{-1} \cap T}) = \mu((g * f) \vert_X) = (\mu \circ E) (g * f)$.
\eproof

Given a trace $\tau$ on $C^*_r(G,\Sigma)$, when do we have $\tau = \tau \circ E$?

\blemma(see also \cite[Proposition~II.5.4]{Ren})
\label{tau=tauE}
Assume that $G$ is principal. Then every trace $\tau$ on $C^*_r(G,\Sigma)$ satisfies $\tau = \tau \circ E$.
\elemma
%\blemma
%\label{tau=tauE}
%Assume that $G$ is principal. Then every trace $\tau$ on $C^*_r(G,\Sigma)$ satisfies $\tau = \tau \circ E$.
%\elemma

\bproof
It suffices to show that for $f \in C_c(G,\Sigma)$ with $f \vert_X \equiv 0$, we have $\tau(f) = 0$. Since $G$ is principal, $G \setminus X$ can be covered by open bisections $S$ with $s(S) \cap r(S) = \emptyset$. So we may assume that $\supp(f) \subseteq S$ for such an open bisection $S$. Choose $h \in C(X)$ such that $\supp(h) \subseteq s(S)$ and $h \vert_{\supp(f)} \equiv 1$. Then $(f * h) (\zeta) = f(\zeta)$ for all $\zeta \in \Sigma$, i.e., $f * h = f$. Moreover, $(h * f) (\zeta) = h(\zeta \eta^{-1}) f(\eta)$ for some $\eta \in j^{-1}(S)$ with $s(j(\eta)) = s(j(\zeta))$. As $s(j(\zeta \eta^{-1})) = s(j(\eta^{-1})) \in s(S^{-1}) = r(S)$ and $h \vert_{r(S)} \equiv 0$, we must have $h(\zeta \eta^{-1}) = 0$ and thus $h * f = 0$. Therefore, $\tau(f) = \tau(f * h) = \tau(h * f) = 0$.
\eproof

\bcor
\label{COR:tau=tauE}
Let $A$ be a separable C*-algebra. Let $B \subseteq A$ be a Cartan subalgebra with faithful conditional expectation $E: \: A \onto B$. Assume that $B$ is a C*-diagonal or that $A$ has unique trace. Then we have $\tau = \tau \circ E$ for every trace $\tau$ on $A$.
\ecor

\subsection{A relationship between C*-diagonals and Cartan subalgebras in vN-algebras}

Let $G$ and $\Tz \times G^{(0)} \overset{i}{\rightarrowtail} \Sigma \overset{j}{\twoheadrightarrow} G$ be as above. Let $c: \: G \to \Sigma$ be a Borel map such that $j \circ c = \id_G$ (i.e., $c$ is a section for $j$), $c \vert_X = \id_X$ and $c(g^{-1}) = c(g)^{-1}$ for all $g \in G$. Existence of such a Borel section was shown in the proof of \cite[Lemma~3.2]{MW}. We then have for every $\zeta \in \Sigma$: $j(\zeta c(j(\zeta))^{-1}) = r(j(\zeta))$, so that there is a Borel map $t: \: \Sigma \to \Tz$ such that $\zeta = t(\zeta) c(j(\zeta))$ for all $\zeta \in \Sigma$. Set $\sigma: \: G^{(2)} \to \Tz, \, (g,h) \ma t(c(gh)^{-1} c(g) c(h))$. $\sigma$ is a normalized cocycle in the sense of \cite[\S~7]{FMI}. Moreover, let $\mu$ be a $\cS(G)$-invariant Borel probability measure on $X$. Set $B \defeq C(X)$, $A \defeq C^*_r(G,\Sigma)$ and let $E: \: A \onto B$ be the unique faithful conditional expectation. Set $\tau \defeq \mu \circ E$. $\tau$ is a trace on $A$ by Lemma~\ref{muE:trace}. Let $\pi \defeq \pi_{\tau}$ be the GNS representation of $A$ attached to $\tau$, and denote by $H \defeq H_{\pi}$ the underlying Hilbert space. Finally, let $R$ be the equivalence relation on $X$ corresponding to $G$, i.e., $R = \menge{(r(g),s(g)) \in X \times X}{g \in G}$.
\blemma
\label{Lem:NoDiag}
$\pi(A)'' \cong \bM(R,\sigma)$.
\elemma
Here $\pi(A)''$ is the von Neumann algebra generated by $\pi(A)$ in $\cL(H)$, and $\bM(R,\sigma)$ is the von Neumann algebra constructed in \cite{FMII}.
\bproof
Let $\nu$ be the right counting measure of $\mu$ as in \cite[Theorem~2]{FMI}. Since we have a Borel isomorphism $G \cong R, \, g \ma (r(g),s(g))$, we may view $\nu$ as a measure on $G$ or $R$, and we have $L^2(R,\nu) \cong L^2(G,\nu)$.
\setlength{\parindent}{0.5cm} \setlength{\parskip}{0cm}

First of all, the map $C_c(G,\Sigma) \to L^2(R,\nu), \, f \ma \dot{f}$, where $\dot{f}(g) = f(c(g))$, extends to a unitary $U: \: H \cong L^2(R,\nu)$. Note that underlying this observation is the isomorphism of Borel groupoids $\Tz \times_{\sigma} G \to \Sigma, \, (z,g) \ma z c(g)$, with inverse given by $(t(\zeta),j(\zeta)) \mafr \zeta$. Here $\Tz \times_{\sigma} G$ is the groupoid with multiplication given by
$$(z_1,\gamma_1)(z_2,\gamma_2) = (z_1 z_2 \sigma(\gamma_1,\gamma_2), \gamma_1 \gamma_2).$$

Secondly, it is straightforward to check that $U \pi(A)'' U^* = \bM(R,\sigma)$.
\eproof

Using \cite[Proposition~5.11]{R08} and Lemma~\ref{tau=tauE}, we obtain
\bcor
Let $A$ be a separable C*-algebra, $B \subseteq A$ a C*-diagonal. Let $\tau$ be a trace on $A$. Let $\pi \defeq \pi_{\tau}$ be the GNS representation of $A$ attached to $\tau$. Then $\pi(A)''$ has a Cartan subalgebra (in the sense of vN-algebras).
\ecor
In other words, if $\pi(A)''$ has no Cartan subalgebra (in the sense of vN-algebras), then $A$ cannot admit C*-diagonals. In particular, this applies to reduced group C*-algebras if the corresponding group vN-algebras have no Cartan subalgebras. Classes of examples of groups with this property are given in \cite{CS,Oza,OP10a,OP10b,PV14a,PV14b}, to mention just a few. For instance, we obtain:
\bcor
\label{Cor:ex-no-Cdiagonal}
Let $\Gamma$ be a discrete group with property (HH)${}^+$ in the sense of \cite{OP10b}, or let $\Gamma$ be an icc hyperbolic group. Then $C^*_r(\Gamma)$ has no C*-diagonal.
\ecor
In particular, this applies to non-abelian free groups. However, we will see later on that for non-abelian free groups, and actually for many more of the groups mentioned in Corollary~\ref{Cor:ex-no-Cdiagonal}, the reduced group C*-algebras have no Cartan subalgebras.

\subsection{Non-existence of Cartan subalgebras in certain reduced group \\ C*-algebras}

Our goal now is to show non-existence of Cartan subalgebras, not only diagonals, in certain reduced group C*-algebras. This builds on recent results about Cartan subalgebras in von Neumann algebras in \cite{BHV}, which go back to \cite{Voi,Hay,Oza,OP10a,OP10b,PV14a,PV14b}.

We start with the following observation: Let $G$ and $\Tz \times G^{(0)} \overset{i}{\rightarrowtail} \Sigma \overset{j}{\twoheadrightarrow} G$ be as above, and let $c: \: G \to \Sigma$ be a Borel section for $j$ with $c \vert_X = \id_X$ and $c(g^{-1}) = c(g)^{-1}$ for all $g \in G$ as before. Assume that $G$ is topologically principal, and let $E: \: C^*_r(G,\Sigma) \onto C(X)$ be the faithful conditional expectation. Let $\tau$ be a trace on $C^*_r(G,\Sigma)$ and $\mu$ the corresponding measure on $X$. Let $(\pi,H_{\pi})$ be the GNS representation of $\tau$. Let $\nu$ be the measure on $G$ given by $\int_G f {\rm d} \nu = \int_X \sum_{\gamma \in G^x} f(\gamma) {\rm d} \mu(x)$. Just as in the proof of Lemma~\ref{Lem:NoDiag}, we obtain
\blemma
\label{H=L^2}
The map $C_c(G,\Sigma) \to C_c(G), \, f \ma f \circ c$ extends to a unitary $H \overset{\cong}{\lori} L^2(G,\nu)$.
\elemma

\blemma
\label{X--L}
We have $\pi(C(X))'' \cong L^{\infty}(X,\mu)$.
\elemma
\bproof
Let $U$ be the unitary $H \overset{\cong}{\lori} L^2(G,\nu)$ from Lemma~\ref{H=L^2}. It is clear that for $f \in C(X)$, $(U \pi(f) U^* \xi) (\gamma) = f(r(\gamma)) \xi(\gamma)$. Since we can cover $G$ by countably many open bisections, we can find a cover of $G$ by countably many pairwise disjoint Borel subsets $S_i \subseteq G$ which are bisections, i.e., $G = \bigsqcup_i S_i$. It is clear that $\nu \vert_{S_i}$ is the pushforward of $\mu \vert_{r(S_i)}$ under $(r \vert_{S_i})^{-1}: \: r(S_i) \overset{\cong}{\lori} S_i$. For $j \in \Nz \cup \gekl{\infty} = \gekl{1, 2, \dotsc} \cup \gekl{\infty}$, define $X_j \defeq \menge{x \in X}{\abs{G^x} = j}$. Clearly, $X_j$ are Borel subsets of $G$ because all the $S_i$s are. Define $I_j \defeq \gekl{1, \dotsc, j}$ for $j \in \Nz$ and $I_j \defeq \Nz$ if $j = \infty$. Then we have a canonical unitary
$$
  V: \: L^2(G,\nu) = \bigoplus_i L^2(S_i,\nu) \cong \bigoplus_i L^2(r(S_i),\mu \vert_{r(S_i)}) \cong \bigoplus_{j \in \Nz \cup \gekl{\infty}} L^2(X_j,\mu \vert_{X_j}) \otimes \ell^2(I_j).
$$
Let $M_j$ be the canonical representation of $L^{\infty}(X_j,\mu \vert_{X_j})$ on $L^2(X_j,\mu \vert_{X_j})$. Then for every $f \in C(X)$, $V \pi(f) V^* = \rukl{M_j(f) \otimes \id_{\ell^2(I_j)}}_{j \in \Nz \cup \gekl{\infty}}$. Hence,
$$
\pi(C(X))'' = \prod_{j \in \Nz \cup \gekl{\infty}} L^{\infty}(X_j,\mu \vert_{X_j}) \cong L^{\infty}(X,\mu).
$$
\eproof

\blemma
\label{LEM:diffuse}
Assume that $\tau$ is a faithful trace on $C^*_r(G,\Sigma)$, and that $\pi(C^*_r(G,\Sigma))''$ is a II$_1$-factor. Then $\mu$ is non-atomic. In particular, $\pi(C(X))''$ is diffuse.
\elemma
\bproof
Assume that $\mu$ is atomic, so that $\mu \gekl{x} > 0$ for some $x \in X$. As $\mu$ is a probability measure on $X$ which is $\cS(G)$-invariant by Lemma~\ref{muSGinv}, we deduce that the orbit $O(x) \defeq \menge{r(\gamma)}{\gamma \in G_x}$ is finite. But that would imply that the characteristic function $1_{O(x)}$ of the orbit of $x$ lies in the center of $\pi(C^*_r(G,\Sigma))''$, which is given by $\Cz 1$ by assumption. Hence $1_{O(x)} = 1$. This means that $X = O(x)$ is finite because $\tau$ is faithful. Hence $G$ must be finite, because otherwise $G$ could not be topologically principal. The conclusion is that $\pi(C^*_r(G,\Sigma))$ is finite dimensional and hence $\pi(C^*_r(G,\Sigma))''$ cannot be a II$_1$-factor.
\eproof

\bcor
\label{Cor:NoCartan}
Let $\Gamma$ be a countable discrete group. Assume that $\Gamma$ is torsion-free. Furthermore, suppose that $\Gamma$ is non-amenable, and that $\Gamma$ has the CMAP and admits a proper 1-cocycle into an orthogonal representation that is weakly contained in the regular representation. Then $C^*_{\lambda}(\Gamma)$ has no Cartan subalgebra.
\ecor
In particular, for $r \geq 2$, $C^*_{\lambda}(\Fz_r)$ has no Cartan subalgebra, and more generally, icc hyperbolic hyperbolic groups have no Cartan subalgebras. More examples where Corollary~\ref{Cor:NoCartan} applies appear in \cite{OP10b}.

\bproof
First of all, our assumptions imply that $C^*_{\lambda}(\Gamma)$ has unique trace: As $\Gamma$ is non-amenable, and because $\Gamma$ has the CMAP and admits a proper 1-cocycle into an orthogonal representation that is weakly contained in the regular representation, we know that $\Gamma$ has the property strong (HH)${}^+$, in the sense of \cite[Definitions~1]{OP10b}. Hence, by \cite[Corollary~2.2]{OP10b}, $\Gamma$ has no infinite normal amenable subgroups. As $\Gamma$ is in addition torsion-free, $\Gamma$ has no non-trivial normal amenable subgroups. In other words, the amenable radical of $\Gamma$ is trivial. By \cite[Theorem~1.3]{BKKO}, we conclude that $C^*_{\lambda}(\Gamma)$ has unique trace.
\setlength{\parindent}{0.5cm} \setlength{\parskip}{0cm}

Write $A = C^*_{\lambda}(\Gamma)$. Assume that $A$ has a Cartan subalgebra $B$, and let $E: \: A \onto B$ be the faithful conditional expectation. Then we can find a groupoid $G$ with compact unit space $X \defeq G^{(0)}$ and a twist $\Sigma$ of $G$ such that $(A,B) \cong (C^*_r(G,\Sigma),C(X))$. Let $\mu$ be the probability measure on $X$ corresponding to the state $\tau \vert_{C(X)}$. Since $C^*_{\lambda}(\Gamma)$ has unique trace, we must have $\tau = \tau \circ E$ (see Corollary~\ref{COR:tau=tauE}). Moreover, let $\pi$ be the GNS representation for $\tau$, so that $\pi(C^*_{\lambda}(\Gamma))'' = L \Gamma$. By Lemma~\ref{X--L}, we have $\pi(B)'' \cong L^{\infty}(X,\mu)$. Since $C^*_{\lambda}(\Gamma)$ has unique trace, $\Gamma$ is icc, so that $L \Gamma$ is a II$_1$-factor. Lemma~\ref{LEM:diffuse} implies that $\pi(B)''$ is diffuse. Hence \cite[Theorem~3.8]{BHV} implies that $\rukl{s \cN_{L \Gamma}(\pi(B)'')}''$ is amenable, where
$$
  s \cN_{L \Gamma}(\pi(B)'') = \menge{x \in L \Gamma}{x \pi(B)'' x^* \subseteq \pi(B)'' \ {\rm and} \ x^* \pi(B)'' x \subseteq \pi(B)''}.
$$
At the same time, it is clear that $s \cN_{L \Gamma}(\pi(B)'')$ contains $\pi(N_A(B))$, where
$$N_A(B) = \menge{n \in A}{nBn^* \subseteq B \ {\rm and} \ n^*Bn \subseteq B}.$$
Since $B$ is a Cartan subalgebra in $A$, we know that $N_A(B)$ generates $A$ as a C*-algebra, so that $L \Gamma = \pi(N_A(B))'' = \rukl{s \cN_{L \Gamma}(\pi(B)'')}''$ is amenable. This cannot be true as $\Gamma$ is not amenable.
\eproof

\bquestion
Is there a C*-algebraic proof for the observation that for $r \geq 2$, $C^*_{\lambda}(\Fz_r)$ has no Cartan subalgebra? 
\equestion

\section{Non-uniqueness of Cartan subalgebras}

Now let us present several classes of C*-algebras which have infinitely many pairwise inequivalent Cartan subalgebras.
\setlength{\parindent}{0cm} \setlength{\parskip}{0.25cm}

The following is a consequence of \cite{DPS}:
\bprop
\label{Prop:ZstableNonUnique}
Let $A$ be a unital separable C*-algebra. Assume that $A$ has a Cartan subalgebra $B$ with $\dim \Spec(B) < \infty$. If $A \cong A \otimes \cZ$, then $A$ has infinitely many pairwise inequivalent Cartan subalgebras.
\eprop
\setlength{\parindent}{0cm} \setlength{\parskip}{0cm}

Here $\dim$ stands for covering dimension, and $\cZ$ is the Jiang-Su algebra.
\bproof
By \cite[Theorem~2.9]{DPS}, the Jiang-Su algebra $\cZ$ has infinitely many Cartan subalgebras $B_1, B_2, \dotsc$ with $\dim \Spec(B_i) < \infty$ for all $i = 1, 2, \dotsc$ and $\lim_{i \to \infty} \dim \Spec B_i = \infty$. Then, by \cite[Lemma~5.1]{BL}, $B \otimes B_1, B \otimes B_2, \dotsc$ are Cartan subalgebras of $A \otimes \cZ$. We have
$$\dim \Spec(B \otimes B_i) = \dim (\Spec(B) \times \Spec(B_i)) \leq \dim \Spec(B) + \dim \Spec(B_i) < \infty$$
for all $i = 1, 2, \dotsc$ Moreover,
$$\dim \Spec(B \otimes B_i) = \dim (\Spec(B) \times \Spec(B_i)) \geq 1 + \dim \Spec(B_i)$$
by \cite[Chapter~2, \S~6.3, (18)]{Fed}. As $\lim_{i \to \infty} 1 + \dim \Spec(B_i) = \infty$, there must be infinitely many $j \in \Nz$ such that $\dim \Spec(B \otimes B_j)$ are pairwise distinct. Hence $A \otimes \cZ$ contains infinitely many non-isomorphic Cartan subalgebras, and thus, so does $A$ as $A \cong A \otimes \cZ$.
\eproof
\bexamples
Here are explicit examples where Proposition~\ref{Prop:ZstableNonUnique} applies:
\begin{itemize}
\item Generalized Bunce-Deddens algebras in the sense of \cite{Orf}, i.e., crossed products $C({\bar \Gamma}) \rtimes_r \Gamma$ attached to odometer actions $\Gamma \curvearrowright {\bar \Gamma}$ of residually finite amenable groups on their profinite completions. $\cZ$-stability follows from \cite{Orf,Phi05} and \cite[Corollary~6.4]{TWW}.
\item Crossed products $C(X) \rtimes_r \Gamma$ attached to free minimal actions of  finitely generated nilpotent infinite groups $\Gamma$ on finite dimensional compact spaces $X$. $\cZ$-stability follows from \cite{SWZ,W12}.
\item Examples of C*-algebras attached to equivalence relations constructed recently in \cite{Put}. By construction, these C*-algebras have Cartan subalgebras isomorphic to algebras of continuous functions on the Cantor space. It is shown in \cite{Put} that these C*-algebras are tracially AF, hence AH algebras of slow dimension growth (by a result of H. Lin, see \cite[Theorem~3.3.5]{Ror}), hence $\cZ$-stable by \cite[Last Corollary]{W12}.
\end{itemize}
\eexamples

\bremark
Since the irrational rotation C*-algebras are $\mathcal Z$-stable, the previous proposition shows that these C*-algebras have infinitely many pairwise inequivalent Cartan subalgebras. A. Kumjian had already shown in his 1980 PhD thesis, published as \cite{Kum84}, that these algebras have diagonals whose  spectrum is a disjoint union of an arbitrary finite number of circles. He gives an explicit construction of these diagonals. The second author thanks him for explaining it to him. First, given a countable subgroup $\Gamma$ of $\Rz$, one defines the groupoid of the action of $\Gamma$ on $\Rz$ by translation
$$G(\Gamma)=\{(x,y)\in\Rz\times\Rz: x-y\in \Gamma\}.$$ 
As a set, it is isomorphic to $\Rz\times\Gamma$; one equips it with the product topology, where $\Gamma$ has the discrete topology. Let $\alpha$ be an irrational number and let $G=G(\Zz+\alpha\Zz)$. Fix a positive element $s\in \Zz+\alpha\Zz$. Since the action of the subgroup $s\Zz$ on $\Rz$ is free and proper, the space $Z(s)=G/G(s\Gamma)$ is a principal $G$-space. Its quotient space is $Z(s)/G=\Rz/s\Zz$. The construction of Section 3 of \cite{R1987} produces a C*-module $E(s)=C^*(Z(s))$ over $C^*(G)$. Its algebra of compact operators $A(s)$ is the C*-algebra of the groupoid $(Z(s)*Z(s))/G$, which is the groupoid of the action of the quotient group $(\Zz+\alpha\Zz)/s\Zz$ on $\Rz/s\Zz$. It is an irrational rotation C*-algebra if $s=p+\alpha q$ with $(p,q)$ relatively prime. Otherwise, it is a matrix algebra over an irrational rotation C*-algebra. The crucial Lemma of Section 9 of \cite{Kum84} constructs an isometry from $E(s+t)$ onto  $E(s)\oplus E(t)$ for positive $s,t\in \Zz+\alpha\Zz$. The existence of such an isometry can also be deduced from the cancellation theorem for projective modules over the irrational rotation C*-algebras of \cite{Rie}. Let $B(s)$ be the canonical diagonal of $A(s)$, viewed as the groupoid C*-algebra $C^*((Z(s)*Z(s))/G)$. Then $A(s+t)$ has a diagonal isomorphic to $B(s)\oplus B(t)$. In particular, $A(1)$, which is the C*-algebra of the rotation $2\pi\alpha$, has a diagonal whose spectrum is the disjoint union of two circles. Repeating the process, one can obtain a diagonal whose spectrum is the disjoint union of an arbitrary finite number of circles. 
\eremark

\bremark
Further examples of C*-algebras with inequivalent Cartan subalgebras can be found in \cite{GPS} (see the second remark after \cite[Theorem~2.4]{GPS}) and \cite{Phi07}. These C*-algebras are given as crossed products attached to certain topological dynamical systems with acting group $\Zz$.
\eremark

We present another non-uniqueness result for Cartan subalgebras. The difference to the previous result is that this time, we will find infinitely many pairwise isomorphic, but inequivalent Cartan subalgebras.
\setlength{\parindent}{0cm} \setlength{\parskip}{0.25cm}

First of all, we need realizations of UCT Kirchberg algebras as groupoid C*-algebras. This has been achieved by various authors. In the following, let us formulate a version which is suitable for our purposes. Moreover, we explain how this version follows from \cite{Sp07a,Sp07b}. The only difference to \cite{Sp07a,Sp07b} is that we want groupoid models (with prescribed additional properties) in the unital case as well.

%\blemma
%\label{Kirchberg-GPD}
%For every unital UCT Kirchberg algebra $A$, there exists an {\'e}tale second countable locally compact Hausdorff groupoid $G$ such that
%\begin{enumerate}
%\item[(i)] $C^*_r(G) \cong A$;
%\item[(ii)] $G^{(0)}$ is homeomorphic to the Cantor space;
%\item[(iii)] for every $x \in G^{(0)}$, the stabilizer group $G^x_x$ is free abelian of rank at most $2$;
%\item[(iv)] the canonical inclusion $C(G^{(0)}) \into C^*_r(G)$ is surjective in $K_0$.
%\end{enumerate}
%\elemma
\setlength{\parindent}{0cm} \setlength{\parskip}{0cm}

\blemma
\label{Kirchberg-GPD}
For every unital UCT Kirchberg algebra $A$, there exists an {\'e}tale second countable locally compact Hausdorff groupoid $G$ such that
\begin{enumerate}
\item[(i)] $C^*_r(G) \cong A$;
\item[(ii)] $G^{(0)}$ is homeomorphic to the Cantor space;
\item[(iii)] $G$ is topologically principal;
\item[(iv)] for every $x \in G^{(0)}$, the stabilizer group $G^x_x$ is free abelian of rank at most $2$.
\end{enumerate}
\elemma
\bproof
This is a slight modification of the construction in \cite{Sp07a}. Set
$$(G_0,g,G_1) \defeq (K_0(A),[1],K_1(A)).$$
Using \cite{Sp07b}, find countable, irreducible graphs $E_0$, $E_1$, $F_0$ and $F_1$ with a (unique) vertex emitting infinitely many edges, together with vertices $v$ in $E_0$ and $w$ in $F_0$ such that
\begin{itemize}
\item $(K_0(\cO(E_0)),[P_v],K_1(\cO(E_0)) \cong (G_0,g,\gekl{0})$;
\item $(K_0(\cO(F_0)),[P_w],K_1(\cO(F_0)) \cong (\Zz,1,\gekl{0})$;
\item $(K_0(\cO(E_1)),K_1(\cO(E_1)) \cong (G_1,\gekl{0})$;
\item $(K_0(\cO(F_1)),K_1(\cO(F_1)) \cong (\gekl{0},\Zz)$.
\end{itemize}
Now apply the construction in \cite{Sp07a}. We obtain a \an{graph-like} (or rather \an{$2$-graph-like}) object $\Omega$ out of the $2$-graphs $E_0 \times F_0$ and $E_1 \times F_1$. For $\Omega$ we have a notion of finite paths, denoted by $\Omega^*$, and a length function $l: \: \Omega^* \to \Nz^2$ (where $\Nz = \gekl{0, 1, 2, \dotsc}$). We also have infinite paths in $\Omega$, and we can equip the set $X$ of infinite paths with a natural topology such that $X$ becomes a totally disconnected, second countable, locally compact Hausdorff space with no isolated points. Moreover, there is an obvious way to concatenate finite paths with infinite paths, denoted by $\mu z$ for a finite path $\mu$ and an infinite path $z$. We can then define the groupoid $\ti G$ consisting of all triples $(x,n,y) \in X \times \Zz^2 \times X$ for which there exist finite paths $\mu$, $\nu$ and an infinite path $z$ such that $x = \mu z$, $y = \nu z$ and $l(\mu) - l(\nu) = n$. Inversion and multiplication are given by $(x,n,y)^{-1} = (y,-n,x)$ and $(x,m,y) (y,n,z) = (x,m+n),z)$. In particular, the unit space of $\ti G$ is given by all triples $(x,0,x)$, $x \in X$, and is canonically homeomorphic to $X$. Spielberg shows in \cite{Sp07a} that ${\ti G}$ is topologically principal and that $C^*_r(\ti G)$ is a stable UCT Kirchberg algebra and that there is canonical inclusion
$$
  i: \: (\cO(E_0) \otimes \cO(F_0)) \oplus (\cO(E_1) \otimes \cO(F_1)) \into C^*_r(\ti G).
$$
We do not need the precise form of $i$, but only the following properties:
\begin{enumerate}
\item[(i$_1$)] $i$ induces an isomorphism in K-theory;
\item[(i$_2$)] for every finite path $\mu \in E_0^*$ in $E_0$, $i(S_{\mu} S_{\mu}^* \otimes P_w)$ is of the form $1_Y$ for some compact open subspace $Y \subseteq X$.
\end{enumerate}
In (i$_2$), $S_{\mu}$ is the partial isometry corresponding to $\mu$ in the graph C*-algebra $\cO(E_0)$, and $1_Y$ is the characteristic function of $Y$, viewed as an element in $C_0(X) \cong C_0({\ti G}^{(0)}) \subseteq C^*_r(\ti G)$.

Let $Z$ be the compact open subset of $X$ such that $1_Z = i(P_v \otimes P_w)$. Then $i$ induces an isomorphism
$$(K_0(C^*_r(\ti G)),[1_Z],K_1(C^*_r(\ti G))) \cong (G_0,g,G_1).$$
Hence the canonical inclusion $C^*_r({\ti G}_Z^Z) \cong 1_Z C^*_r(\ti G) 1_Z \into C^*_r(G)$ induces an isomorphism 
$$(K_0(C^*_r({\ti G}_Z^Z)),[1],K_1(C^*_r({\ti G}_Z^Z))) \cong (G_0,g,G_1).$$
Therefore, we may set $G \defeq {\ti G}_Z^Z$, and the classification theorem due to Kirchberg and Phillips \cite[Theorem~8.4.1]{Ror} yields $C^*_r(G) \cong A$. The unit space of $G$ is obviously given by $Z$ and hence homeomorphic to the Cantor space as $Z$ is totally disconnected, second countable, compact Hausdorff and has no isolated points. 
%Moreover, $G$ has property (iii) because $\ti G$ has the same property. 
Moreover, $G$ has properties (iii) and (iv) because it is equivalent to $\ti G$ which has the same properties.
\eproof

With the help of this lemma, we obtain the desired non-uniqueness result for Cartan subalgebras in unital UCT Kirchberg algebras.
\bprop
\label{unitalKirchberg-nonunique}
Every unital UCT Kirchberg algebra has infinitely many pairwise inequivalent Cartan subalgebras, whose spectra are all homeomorphic to the Cantor space.
\eprop
\bproof
Let $A$ be a unital UCT Kirchberg algebra. By Lemma~\ref{Kirchberg-GPD}, we can find a groupoid $G$ with $G^{(0)}$ homeomorphic to the Cantor space such that $C^*_r(G) \cong A$ and for every $x \in G^{(0)}$, the stabilizer group $G_x^x = \menge{\gamma \in G}{r(\gamma) = s(\gamma) = x}$ is a free abelian group with $\sup_{x \in G^{(0)}} \rk G_x^x \leq 2$. We set $\rk_{\rm Stab}(G) \defeq \sup_{x \in G^{(0)}} \rk G_x^x$. In particular, for $\cO_{\infty}$, we obtain a Cartan subalgebra $B_{\infty}$ and a groupoid $G_{\infty}$ such that $(\cO_{\infty},B_{\infty}) \cong (C^*_r(G_{\infty}),C(G_{\infty}^{(0)}))$. We can choose $G_{\infty}$ such that $G_{\infty}^{(0)}$ is homeomorphic to the Cantor space and $\rk_{\rm Stab}(G_{\infty}) = 1$ (see for instance \cite[Chapter~III, \S~2]{Ren}).
\setlength{\parindent}{0.5cm} \setlength{\parskip}{0cm}

For every $i = 1, 2, \dotsc$, we know that $A \cong A \otimes \cO_{\infty}^{\otimes i}$ because $A$ is purely infinite. By \cite[Lemma~5.1]{BL}, we have $A \otimes \cO_{\infty}^{\otimes i} \cong C^*_r(G \times G_{\infty}^i)$, where $G_{\infty}^i$ means here the $i$-fold Cartesian product. In this way, we obtain Cartan subalgebras $B_i$ of $A$ corresponding to $C(G^{(0)} \times (G_{\infty}^{(0)})^i)$ under the isomorphism $A \cong A \otimes \cO_{\infty}^{\otimes i} \cong C^*_r(G \times G_{\infty}^i)$. Obviously, $\Spec(B_i)$ is homeomorphic to the Cantor space for all $i = 1, 2, \dotsc$ As
$$\rk_{\rm Stab}(G \times G_{\infty}^i) = \rk_{\rm Stab}(G) + i \cdot \rk_{\rm Stab}(G_{\infty}) = \rk_{\rm Stab}(G) + i,$$
we know that the groupoids $\menge{G \times G_{\infty}^i}{i = 1, 2, \dotsc}$ must be pairwise non-isomorphic. Hence $(A,B_i)$, $i = 1, 2, \dotsc$, must be pairwise inequivalent.
\eproof

An analogous non-uniqueness result for Cartan subalgebras holds for stable UCT Kirchberg algebras. The proof is the same as for Proposition~\ref{unitalKirchberg-nonunique}. It is even simpler as we can directly use the groupoid models from \cite{Sp07a,Sp07b}.

%\bprop
%Every unital UCT Kirchberg algebra has infinitely many pairwise inequivalent Cartan subalgebras, whose spectra are all homeomorphic to the (up to homeomorphism) unique totally disconnected, second countable, locally compact Hausdorff space without isolated points.
%\eprop

\bprop
Every stable UCT Kirchberg algebra has infinitely many pairwise inequivalent Cartan subalgebras, whose spectra are all homeomorphic to the (up to homeomorphism) unique totally disconnected, second countable, locally compact Hausdorff space without isolated points.
\eprop

\bremark
In \cite{BCSS}, many (but not all) UCT Kirchberg algebras have been realized as C*-algebras of principal -- not only topologically principal -- groupoids.
\eremark

\bquestion
Does every classifiable C*-algebra (i.e., separable, unital, simple, with finite nuclear dimension) have a Cartan subalgebra? 
\equestion

\section{Cartan subalgebras in Roe algebras}

Let us now observe that nuclear Roe algebras have distinguished Cartan subalgebras (see Theorem~\ref{Roe-distCartan}). This follows from \cite{SW} and \cite[\S~2.5]{Li_DQH} (see also \cite{MST}). We refer to \cite{SW,Li_DQH} and the references therein for background material concerning notions from geometric group theory such as quasi-isometry or bilipschitz equivalence.
\setlength{\parindent}{0cm} \setlength{\parskip}{0.25cm}

%\blue{
%Is the following remark OK? Do we want to give more details?
%}
%
%\red{Second countability is used nowhere in Section 4 (except in the implication in Theorem 4.2.(ii) which says that $B$ masa implies $G$ topologically principal and which is not relevant to our purpose) until Proposition 4.13. In the proof of Proposition 4.13, the existence of $h\in C_0(G^{(0)})$ such that the open subset $U=s(S)$ is exactly the subset $h(x)\not=0$ is guaranteed by the $\sigma$-compactness of $G$ (then $G^{(0)}$ and $U$ are also $\sigma$-compact). The density of $C_c(U)$ in $C_0(U)$ in the uniform norm holds for any locally compact Hausdorff space $U$. The theorem of Douady and Soglio-H\'erault is valid under the sole hypothesis that the base space of the Banach bundle is either paracompact or locally compact (in our case, we have both).  Therefore, the proof is valid without any change if $G$ is $\sigma$-compact. Thus, we can simply state that the proof of \cite[Proposition~4.13]{R08} works without any change for \'etale Hausdorff locally compact topologically principal groupoids which are $\sigma$-compact. If you think it is useful, we can mention that $\sigma$-compactness is used to get the existence of the function $h$ in the proof).}

Let us start with the remark that although in \cite{R08}, the groupoids are assumed to be second countable, the proof of \cite[Proposition~4.13]{R08} works without any change for \'etale Hausdorff locally compact topologically principal groupoids which are $\sigma$-compact. Hence, combining this observation with \cite[\S~2.5]{Li_DQH}, we obtain the following result. Let $\cK$ be the C*-algebra of compact operators on $\ell^2 \Zz$ and $\cD \defeq \cK \cap \ell^{\infty}(\Zz)$.
\blemma
\label{qibilip}
\setlength{\parindent}{0cm} \setlength{\parskip}{0cm}

Let $\Gamma$, $\Lambda$ be finitely generated groups. The following are equivalent:
\begin{enumerate}
\item[(QI1)] $\Gamma$ and $\Lambda$ are quasi-isometric;
\item[(QI2)] $\Gamma \curvearrowright \beta \Gamma$ and $\Lambda \curvearrowright \beta \Lambda$ are Kakutani equivalent, in the sense of \cite{Mat,Li_DQH};
\item[(QI3)] $\Gamma \curvearrowright \beta \Gamma$ and $\Lambda \curvearrowright \beta \Lambda$ are stably continuously orbit equivalent, in the sense of \cite{CRS,Li_DQH};
\item[(QI4)] There exist projections $p \in C(\beta \Gamma) \subseteq C(\beta \Gamma) \rtimes_r \Gamma$ and $q \in C(\beta \Lambda) \subseteq C(\beta \Lambda) \rtimes_r \Lambda$ which are full in $C(\beta \Gamma) \rtimes_r \Gamma$ and $C(\beta \Lambda) \rtimes_r \Lambda$, respectively, such that
$$(p(C(\beta \Gamma) \rtimes_r \Gamma)p, pC(\beta \Gamma)p) \cong (q(C(\beta \Lambda) \rtimes_r \Lambda)q, qC(\beta \Lambda)q);$$
\item[(QI5)] $((C(\beta \Gamma) \rtimes_r \Gamma) \otimes \cK, C(\beta \Gamma) \otimes \cD) \cong ((C(\beta \Lambda) \rtimes_r \Lambda) \otimes \cK, C(\beta \Lambda) \otimes \cD)$.
\end{enumerate}

Moreover, the following are equivalent:
\begin{enumerate}
\item[(LIP1)] $\Gamma$ and $\Lambda$ are bilipschitz equivalent;
\item[(LIP2)] $\Gamma \curvearrowright \beta \Gamma$ and $\Lambda \curvearrowright \beta \Lambda$ are continuously orbit equivalent, in the sense of \cite{Li_COER};
\item[(LIP3)] $(C(\beta \Gamma) \rtimes_r \Gamma, C(\beta \Gamma)) \cong (C(\beta \Lambda) \rtimes_r \Lambda, C(\beta \Lambda))$.
\end{enumerate}
\elemma
\setlength{\parindent}{0cm} \setlength{\parskip}{0cm}

Here, $\Gamma \curvearrowright \beta \Gamma$ and $\Lambda \curvearrowright \beta \Lambda$ are the canonical actions of $\Gamma$ and $\Lambda$ on their Stone-{\v C}ech compactifications.
\setlength{\parindent}{0cm} \setlength{\parskip}{0.25cm}

Let us now combine Lemma~\ref{qibilip} with \cite{SW}. First of all, recall that for a finitely generated group $\Gamma$, the uniform Roe algebra $C^*_u(X_{\Gamma})$ is canonically isomorphic to $C(\beta \Gamma) \rtimes_r \Gamma$, and that the stable uniform Roe algebra $C^*_s(X_{\Gamma})$ is canonically isomorphic to $(C(\beta \Gamma) \rtimes_r \Gamma) \otimes \cK$. We let $D_u(X_{\Gamma})$ be the canonical Cartan subalgebra of $C^*_u(X_{\Gamma})$ corresponding to the Cartan subalgebra $C(\beta \Gamma)$ of $C(\beta \Gamma) \rtimes_r \Gamma$. Similarly, $D_s(X_{\Gamma})$ denotes the canonical Cartan subalgebra of $C^*_s(X_{\Gamma})$ corresponding to the Cartan subalgebra $C(\beta \Gamma)\otimes \cD$ of $(C(\beta \Gamma) \rtimes_r \Gamma) \otimes \cK$.
\setlength{\parindent}{0cm} \setlength{\parskip}{0cm}

\btheo
\label{Roe-distCartan}
Let $\Gamma$, $\Lambda$ be finitely generated groups. Consider the following statements:
\begin{enumerate}
\item[(i)] $\Gamma$ and $\Lambda$ are bilipschitz equivalent;
\item[(ii)] $(C^*_u(X_{\Gamma}), D_u(X_{\Gamma})) \cong (C^*_u(X_{\Lambda}), D_u(X_{\Lambda}))$;
\item[(iii)] $C^*_u(X_{\Gamma}) \cong C^*_u(X_{\Lambda})$;
\item[(iv)] $\Gamma$ and $\Lambda$ are quasi-isometric;
\item[(v)] There exist projections $p \in D_u(X_{\Gamma}) \subseteq C^*_u(X_{\Gamma})$ and $q \in D_u(X_{\Lambda}) \subseteq C^*_u(X_{\Lambda})$ which are full in $C^*_u(X_{\Gamma})$ and $C^*_u(X_{\Lambda})$, respectively, such that
$$(pC^*_u(X_{\Gamma})p, pD_u(X_{\Gamma})p) \cong (qC^*_u(X_{\Lambda})q, qD_u(X_{\Lambda})q);$$
\item[(vi)] $(C^*_s(X_{\Gamma}), D_s(X_{\Gamma})) \cong (C^*_s(X_{\Lambda}), D_s(X_{\Lambda}))$;
\item[(vii)] $C^*_s(X_{\Gamma}) \cong C^*_s(X_{\Lambda})$.
\end{enumerate}
We always have (i) $\LRarr$ (ii) and (iv) $\LRarr$ (v) $\LRarr$ (vi). Moreover, (ii) $\Rarr$ (iii), (vi) $\Rarr$ (vii), (ii) $\Rarr$ (v) and (iii) $\Rarr$ (vii) are obvious. If $\Gamma$ and $\Lambda$ are exact, we have (vii) $\Rarr$ (iv), so that (iv) $\LRarr$ (v) $\LRarr$ (vi) $\LRarr$ (vii). For non-amenable groups $\Gamma$ and $\Lambda$, we have (iv) $\Rarr$ (i). Hence for non-amenable, exact groups $\Gamma$ and $\Lambda$, all these items (i) -- (vii) are equivalent.
\etheo
\bproof
(i) $\LRarr$ (ii) and (iv) $\LRarr$ (v) $\LRarr$ (vi) follow from Lemma~\ref{qibilip}. For exact groups, (vii) $\Rarr$ (iv) is \cite[Corollary~6.3]{SW}. Moreover, \cite{Why} tells us that (iv) $\Rarr$ (i) holds for non-amenable groups.
\eproof

\bremark
We know by \cite{Dym} that for amenable groups, (iv) does not imply (i) in general. It would be interesting to find examples of groups where (v) holds but not (iii), or where (iii) holds but not (ii). We know that for the groups in \cite{Dym}, at least one of the implications (v) $\Rarr$ (iii) or (iii) $\Rarr$ (ii) must fail. But which one(s)?
\eremark


\begin{thebibliography}{99}

\bibitem{BL} S. \textsc{Barlak} and X. \textsc{Li}, \emph{Cartan subalgebras and the UCT problem}, preprint, arXiv:1511.02697v2.

\bibitem{Bla} B. \textsc{Blackadar},
	\emph{Operator Algebras, Theory of C*-Algebras and von Neumann Algebras},
	Encyclopaedia of Mathematical Sciences, Vol. 122,
	Springer-Verlag, Berlin Heidelberg, 2006.

\bibitem{BHV} R. \textsc{Boutonnet}, C. \textsc{Houdayer} and S. \textsc{Vaes}, \emph{Strong solidity of free Araki-Woods factors}, arXiv:1512.04820.

\bibitem{BM} R. \textsc{Boyer} and R. \textsc{Martin}, \emph{The regular group C*-algebra for real-rank one groups}, Proc. Amer. Math. Soc. 59 (1976), 371--376. 

\bibitem{BKKO} E. \textsc{Breuillard}, M. \textsc{Kalantar}, M. \textsc{Kennedy} and N. \textsc{Ozawa}, \emph{C*-simplicity and the unique trace property for discrete groups}, preprint, arXiv:1410.2518v3.

\bibitem{BCSS} J.H. \textsc{Brown}, L.O. \textsc{Clark}, A. \textsc{Sierakowski} and A. \textsc{Sims}, \emph{Purely infinite simple C*-algebras that are principal groupoid C*-algebras}, J. Math. Anal. Appl. \emph{439} (2016), 213--234. 

\bibitem{CRS} T.M. \textsc{Carlsen}, E. \textsc{Ruiz} and A. \textsc{Sims}, \emph{Equivalence and stable isomorphism of groupoids, and diagonal-preserving stable isomorphisms of graph $C^*$-algebras and Leavitt path algebras}, Proc. Amer. Math. Soc. \emph{145} (2017), 1581--1592.

\bibitem{CS} I. \textsc{Chifan} and T. \textsc{Sinclair}, \emph{On the structural theory of ${\rm II}_1$ factors of negatively curved groups}, Annales scientifiques de l'{\'E}.N.S., S{\'e}r. 4 \emph{46} (2013), p. 1--33.

\bibitem{DPS} R.J. \textsc{Deeley}, I.F. \textsc{Putnam}, K.R. \textsc{Strung}, \emph{Constructing minimal homeomorphisms on point-like spaces and a dynamical presentation of the Jiang-Su algebra}, preprint, arXiv:1503.03800v2, to appear in J. reine angew. Math.

\bibitem{Dix} J. \textsc{Dixmier}, \emph{C*-algebras}, Translated from the French by Francis Jellett, North-Holland Mathematical Library, Vol. 15, North-Holland Publishing Co., Amsterdam-New York-Oxford, 1977.

\bibitem{DF1} R.S. \textsc{Doran} and J.M.G. \textsc{Fell}, \emph{Representations of *-algebras, locally compact groups, and Banach *-algebraic bundles}, Vol. 1, Basic representation theory of groups and algebras. Pure and Applied Mathematics, 125. Academic Press, Inc., Boston, MA, 1988.

\bibitem{DF2} R.S. \textsc{Doran} and J.M.G. \textsc{Fell}, \emph{Representations of *-algebras, locally compact groups, and Banach *-algebraic bundles}, Vol. 2, Banach *-algebraic bundles, induced representations, and the generalized Mackey analysis. Pure and Applied Mathematics, 126. Academic Press, Inc., Boston, MA, 1988.

\bibitem{Dym} T. \textsc{Dymarz}, \emph{Bilipschitz equivalence is not equivalent to quasi-isometric equivalence for finitely generated groups}, Duke Math. J. \emph{154} (2010), 509--526.

\bibitem{EGLN} G. \textsc{Elliott}, G. \textsc{Gong}, H. \textsc{Lin} and Z. \textsc{Niu}, \emph{On the classification of simple amenable C*-algebras with finite decomposition rank, II}, preprint, arXiv:1507.03437v3.

\bibitem{Fed} V.V. \textsc{Fedorchuk}, \emph{The Fundamentals of Dimension Theory}, in Encyclopaedia of Mathematical Sciences, Vol. 17, General Topology I, A.V. Arkhangel'skii and L.S. Pontryagin (Eds.), Springer-Verlag, Berlin, 1993.

\bibitem{FMI} J. \textsc{Feldman} and C.C. \textsc{Moore}, \emph{Ergodic Equivalence Relations, Cohomology, and von Neumann Algebras. I}, Trans. Amer. Math. Soc. \emph{234} (1977), 289--324.

\bibitem{FMII} J. \textsc{Feldman} and C.C. \textsc{Moore}, \emph{Ergodic Equivalence Relations, Cohomology, and von Neumann Algebras. II}, Trans. Amer. Math. Soc. \emph{234} (1977), 325--359.

\bibitem{Fel} J.M.G. \textsc{Fell}, \emph{The structure of algebras of operator fields}, Acta Math. \emph{106} (1961), 233--280.

\bibitem{Fur11} A. \textsc{Furman}, \emph{A survey of measured group theory}, in Geometry, rigidity, and group actions, 296--374, Chicago Lectures in Math., Univ. Chicago Press, Chicago, IL, 2011.

\bibitem{Gab10} D. \textsc{Gaboriau}, \emph{Orbit equivalence and measured group theory}, Proc. Int. Congr. Math., Volume III, 1501--1527, Hindustan Book Agency, New Delhi, 2010.

\bibitem{GPS} T. \textsc{Giordano}, I.F. \textsc{Putnam} and C.F. Skau, \emph{Topological orbit equivalence and $C^*$-crossed products}, J. Reine Angew. Math. \emph{469} (1995), 51--111.

\bibitem{GLN} G. \textsc{Gong}, H. \textsc{Lin} and Z. \textsc{Niu}, \emph{Classification of finite simple amenable $\cZ$-stable C*-algebras}, preprint, arXiv:1501.00135v6.

\bibitem{Gre} K.D. \textsc{Gregson}, \emph{Extension of Pure States of C*-Algebras}, Ph.D. thesis, University of Aberdeen, Aberdeen, 1986.

\bibitem{Hay} B. \textsc{Hayes}, \emph{1-bounded entropy and regularity problems in von Neumann algebras}, to appear in Int. Math. Res. Not., arXiv:1505.06682v4.

\bibitem{Hus} D. \textsc{Husemoller}, \emph{Fibre bundles}, Third edition, Graduate Texts in Mathematics, 20, Springer-Verlag, New York, 1994.

\bibitem{KM} H. \textsc{Kraljevi{\'c}} and D. \textsc{Mili{\v c}i{\'c}}, \emph{The C*-algebra of the universal covering group of SL(2,R)}, Glasnik Mat. Ser. III \emph{7} (1972), 35--48.

\bibitem{KL} F. \textsc{Krauss} and T.C. \textsc{Lawson}, \emph{Examples of homogeneous C*-algebras}, Mem. Amer. Math. Soc. \emph{148} (1974), 153--164.

%added:

%--------------------------------

\bibitem{Kum84} A. \textsc{Kumjian}, \emph{On localizations and simple C*-algebras}, Pacific J. Math.\emph{112} (1984), 141--192.

%------------------------------------

\bibitem{Kum} A. \textsc{Kumjian}, \emph{Diagonals in algebras of continuous trace. With an appendix by Toshikazu Natsume}, Lecture Notes in Math., 1132, Operator algebras and their connections with topology and ergodic theory (Bu{\c s}teni, 1983), 297--311, Springer, Berlin, 1985.

\bibitem{Li_COER} X. \textsc{Li}, \emph{Continuous orbit equivalence rigidity}, to appear in Ergod. Th. Dyn. Sys., arXiv:1503.01704.

\bibitem{Li_PTGPD} X. \textsc{Li}, \emph{Partial transformation groupoids attached to graphs and semigroups}, preprint, arXiv: 1603.09165, accepted for publication in Int. Math. Res. Not.

\bibitem{Li_DQH} X. \textsc{Li}, \emph{Dynamic characterizations of quasi-isometry, and applications to cohomology}, preprint, arXiv: 1604.07375v3.

\bibitem{Mat} H. \textsc{Matui}, \emph{Homology and topological full groups of {\'e}tale groupoids on totally disconnected spaces}, Proc. Lond. Math. Soc. \emph{104} (2012), 27--56.

\bibitem{MS12} H. \textsc{Matui} and Y. \textsc{Sato}, \emph{Strict comparison and $\cZ$-absorption of nuclear C*-algebras}, Acta Math. \emph{209} (2012), 179--196. 

\bibitem{MS14} H. \textsc{Matui} and Y. \textsc{Sato}, \emph{Decomposition rank of UHF-absorbing C*-algebras}, Duke Math. J. \emph{163} (2014), 2687--2708.

\bibitem{MST} K. \textsc{Medynets}, R. \textsc{Sauer} and A. \textsc{Thom}, \emph{Cantor systems and quasi-isometry of groups}, preprint, arXiv:1508.07578v2.

\bibitem{Mil} D. \textsc{Mili{\v c}i{\'c}}, \emph{Topological representation of the group C*-algebra of SL(2,R)}, Glasnik Mat. Ser. III \emph{6} (1971), 231--246.

\bibitem{MW} P.S. \textsc{Muhly} and D.P. \textsc{Williams}, \emph{Continuous trace groupoid C*-algebras, II}, Math. Scand. \emph{70} (1992), 127--145.

\bibitem{Orf} S. \textsc{Orfanos}, \emph{Generalized Bunce-Deddens algebras}, Proc. Amer. Math. Soc. \emph{138} (2010), 299--308.

\bibitem{Oza} N. \textsc{Ozawa}, \emph{Solid von Neumann algebras}, Acta Math. \emph{192} (2004), 111--117.

\bibitem{OP10a} N. \textsc{Ozawa} and S. \textsc{Popa}, \emph{On a class of ${\rm II}_1$ factors with at most one Cartan subalgebra}, Ann. of Math. \emph{172} (2010), 713--749.

\bibitem{OP10b} N. \textsc{Ozawa} and S. \textsc{Popa}, \emph{On a class of II1 factors with at most one Cartan subalgebra, II}, Amer. J. Math. \emph{132} (2010), 841--866.

\bibitem{PP} M.B. \textsc{Penington} and R. \textsc{Plymen}, \emph{The Dirac operator and the principal series for complex semisimple Lie groups}, J. Funct. Anal. \emph{53} (1983), 269--286. 

\bibitem{Phi05} N.C. \textsc{Phillips}, \emph{Crossed products of the Cantor set by free minimal actions of $\Zz^d$}, Comm. Math. Phys. \emph{256} (2005), 1--42. 

\bibitem{Phi07} N.C. \textsc{Phillips}, \emph{Examples of different minimal diffeomorphisms giving the same $C^*$-algebras}, Israel J. Math. \emph{160} (2007), 189--217. 

\bibitem{PV14a} S. \textsc{Popa} and S. \textsc{Vaes}, \emph{Unique Cartan decomposition for II1 factors arising from arbitrary actions of free groups}, Acta Math. \emph{212} (2014), 141--198.

\bibitem{PV14b} S. \textsc{Popa} and S. \textsc{Vaes}, \emph{Unique Cartan decomposition for II${}_1$ factors arising from arbitrary actions of hyperbolic groups}, J. reine angewandte Math. \emph{694} (2014), 215--239.

\bibitem{Put} I. \textsc{Putnam}, \emph{Some classifiable groupoid C*-algebras with prescribed K-theory}, preprint, arXiv: 1611.04649v2.

\bibitem{Ren} J. \textsc{Renault}, \emph{A groupoid approach to C*-algebras}, Lecture Notes in Math., 793, Springer, Berlin, 1980.
 
\bibitem{R1983} J. \textsc{Renault}, \emph{Two applications of the dual groupoid of a C*-algebra}, Lecture Notes in Math., 1132, Operator algebras and their connections with topology and ergodic theory (Bu{\c s}teni, 1983), 434--445, Springer, Berlin, 1985. 

%added:

%--------------------------------

\bibitem{R1987} J. \textsc{Renault}, \emph{Repr\'esentations des produits crois\'es d'alg\`ebres de groupo\"ides}, J. Oper. Th. \emph{18} (1987), 67--97.

%------------------------------------


\bibitem{R08} J. \textsc{Renault}, \emph{Cartan subalgebras in C*-algebras}, Irish Math. Soc. Bull. \emph{61} (2008), 29--63.

%added:

%--------------------------------

\bibitem{Rie} M. \textsc{Rieffel}, \emph{The cancellation theorem for projective modules over the irrational rotation C*-algebras}, Proc. London Math. Soc. (3) \emph{47} (1983), no.2, 285--302.

%------------------------------------



\bibitem{Ror} M. \textsc{R{\o}rdam},
	\emph{Classification of Nuclear $C^*$-Algebras}
	in Classification of Nuclear $C^*$-Algebras. Entropy in Operator Algebras,
	Encyclopaedia of Mathematical Sciences, Vol. 126, Springer-Verlag,
	Berlin Heidelberg New York, 2002.

\bibitem{SWW} Y. \textsc{Sato}, S. \textsc{White} and W. \textsc{Winter}, \emph{Nuclear dimension and $\cZ$-stability}, Invent. Math. \emph{202} (2015), 893--921.

\bibitem{Sha05} Y. \textsc{Shalom}, \emph{Measurable group theory}, Eur. Congr. Math., 391--423, Eur. Math. Soc., Z{\"u}rich, 2005.

\bibitem{SW} J. \textsc{Spakula} and R. \textsc{Willett}, \emph{On rigidity of Roe algebras}, Adv. Math. \emph{249} (2013), 289--310.

\bibitem{Sp07a} J. \textsc{Spielberg}, \emph{Graph-based models for Kirchberg algebras}, J. Operator Th. \emph{57} (2007), 347--374. 

\bibitem{Sp07b} J. \textsc{Spielberg}, \emph{Non-cyclotomic presentations of modules and prime-order automorphisms of Kirchberg algebras}, J. Reine Angew. Math. \emph{613} (2007), 211--230. 

\bibitem{Ste} N. \textsc{Steenrod}, \emph{The topology of fibre bundles}, Landmarks in Mathematics, Princeton Paperbacks, Princeton University Press, Princeton, NJ, 1951.

\bibitem{SWZ} G. \textsc{Szabo}, J. \textsc{Wu} and J. \text{Zacharias}, \emph{Rokhlin dimension for actions of residually finite groups}, preprint, arXiv:1408.6096v4.

\bibitem{TWW} A. \textsc{Tikuisis}, S. \textsc{White} and W. \textsc{Winter}, \emph{Quasidiagonality of nuclear C*-algebras},  Ann. of Math. (2) \emph{185} (2017), 229--284.

\bibitem{Val} A. \textsc{Valette}, \emph{K-theory for the reduced C*-algebra of a semisimple Lie group with real rank 1 and finite centre}, Quart. J. Math. Oxford Ser. (2) \emph{35} (1984), 341--359. 

\bibitem{Voi} D. \textsc{Voiculescu}, \emph{The analogues of entropy and of Fisher's information measure in free probability theory. III. The absence of Cartan subalgebras}, Geom. Funct. Anal. \emph{6} (1996), 172--199. 

\bibitem{Why} K. \textsc{Whyte}, \emph{Amenability, bilipschitz equivalence, and the von Neumann conjecture}, Duke Math. J. \emph{99} (1999), 93--112.

\bibitem{W10} W. \textsc{Winter}, \emph{Decomposition rank and $\cZ$-stability}, Invent. Math. \emph{179} (2010), 229--301.

\bibitem{W12} W. \textsc{Winter}, \emph{Nuclear dimension and $\cZ$-stability of pure C*-algebras}, Invent. Math. \emph{187} (2012), 259--342.

\end{thebibliography}
\end{document}